\documentclass[final,onefignum,onetabnum]{siamonline190516}

\usepackage{lipsum}
\usepackage{amsfonts}
\usepackage{epstopdf}
\ifpdf
  \DeclareGraphicsExtensions{.eps,.pdf,.png,.jpg}
\else
  \DeclareGraphicsExtensions{.eps}
\fi

\usepackage{datetime}
\newdateformat{monthyeardate}{%
  \monthname[\THEMONTH], \THEYEAR}

\usepackage{amsmath}
\allowdisplaybreaks
\usepackage{amssymb}
\usepackage{commath}
\usepackage{mathtools}
\usepackage{bbm}

\usepackage{color}
\usepackage{graphicx}
\usepackage[small]{caption}
\usepackage{subcaption}

\usepackage{relsize}
\usepackage{adjustbox}
\usepackage{algorithm}
\usepackage[noend]{algpseudocode}
\usepackage{booktabs}
\usepackage{tikz}

\usepackage{enumitem}
\setlist[enumerate]{leftmargin=.5in}
\setlist[itemize]{leftmargin=.5in}

\newsiamremark{remark}{Remark}
\newsiamremark{hypothesis}{Hypothesis}
\crefname{hypothesis}{Hypothesis}{Hypotheses}
\newsiamthm{claim}{Claim}

\headers{Constructing Positive Interpolatory CFs}{Jan Glaubitz}

\title{Constructing Positive Interpolatory Cubature Formulas\thanks{\monthyeardate\today \funding{This work was supported by the German Research Foundation (DFG) under grant GL 927/1-1.}}}

\author{Jan Glaubitz\thanks{Department of Mathematics, Dartmouth College, Hanover, NH 03755, USA 
  (\email{Jan.Glaubitz@Dartmouth.edu}, \url{https://math.dartmouth.edu/\~jglaubitz}).}
}

\usepackage{amsopn}

\DeclareMathOperator{\diag}{diag}
\DeclareMathOperator*{\argmin}{arg\,min} 
\newcommand{\scp}[2]{\left\langle{#1,\, #2}\right\rangle} 
\newcommand{\intd}{\, \mathrm{d}}
\newcommand{\N}{\mathbb{N}}

\newcommand{\R}{\mathbb{R}}

\ifpdf
\hypersetup{
  pdftitle={Constructing Positive Interpolatory Cubature Formulas},
  pdfauthor={Jan Glaubitz}
}
\fi

\begin{document}

\maketitle

% REQUIRED
\begin{abstract}
	Positive interpolatory cubature formulas (CFs) are constructed for quite general integration domains and weight functions. 
These CFs are exact for general vector spaces of continuous real-valued functions that contain constants. 
At the same time, the number of data points---all of which lie inside the domain of integration--- and cubature weights---all positive---is less or equal to the dimension of that vector space. 
The existence of such CFs has been ensured by Tchakaloff in 1957. 
Yet, to the best of the author's knowledge, this work is the first to provide a procedure to successfully construct them. 
\end{abstract}

% REQUIRED

% REQUIRED

\section{Introduction} 
\label{sec:introduction} 

% Introduction 
Dating back as far as to the ancient Babylonians and Egyptians \cite{boyer2011history}, numerical integration has always been an omnipresent problem in mathematics. 
Let $d \geq 2$, and $\Omega \subset \R^d$ be a compact set with positive volume $|\Omega|$. 
Then, the task is to approximate the weighted integral 
\begin{equation}\label{eq:int}
	I[f] := \int_\Omega \omega(\boldsymbol{x}) f(\boldsymbol{x}) \intd \boldsymbol{x}
\end{equation} 
with nonnegative \emph{weight function} $\omega: \Omega \to \R_0^+$ by an \emph{$N$-point CF}
\begin{equation}\label{eq:CF}
	C_N[f] := \sum_{n=1}^N w_n f(\mathbf{x}_n). 
\end{equation} 
The distinct points $\{ \mathbf{x}_n \}_{n=1}^N \subset \R^d$ are called \emph{data points} and the coefficients $\{ w_n \}_{n=1}^N \subset \R$ are referred to as \emph{cubature weights}. 
% Points inside and weights nonnegative
For a 'good' CF, the following properties are often required:
\begin{enumerate}[label=(P\arabic*)] 
	\item \label{item:P1}
	All data points should lie inside of $\Omega$. That is, $\mathbf{x}_n \in \Omega$ for all $n=1,\dots,N$. 
	
	\item \label{item:P2} 
	All cubature weights should be positive. That is, $w_n > 0$ for all $n=1,\dots,N$.\footnote{Some authors require the cubature weights to only be nonnegative. 
Yet, for $w_n = 0$, the corresponding data point can and should be removed from the CF to avoid an unnecessary loss of efficiency.}
\end{enumerate}  
See the many excellent monographs and reviews \cite{haber1970numerical,stroud1971approximate,mysovskikh1980approximation,engels1980numerical,cools1997constructing,davis2007methods} on CFs.
% Exactness 
Furthermore, it is often desired that a CF is exact for certain linear vector spaces. 
Let $\mathcal{F}_K(\Omega)$ be the $K$-dimensional linear vector space spanned by the linearly independent continuous functions ${\varphi_1,\dots,\varphi_K:\Omega \to \R}$. 
It shall be assumed that all the moments $m_k := I[\varphi_k]$, $k=1,\dots,K$, exist.
Then, the $N$-point CF \eqref{eq:CF} is said to be \emph{$\mathcal{F}_K(\Omega)$-exact} if 
\begin{equation}
	C_N[f] = I[f] \quad \forall f \in \mathcal{F}_K(\Omega). 
\end{equation} 
Usual choices for $\mathcal{F}_K(\Omega)$ include the vector space of all algebraic and trigonometric polynomials up to a certain degree $m$, respectively denoted by $\mathbb{P}_m(\R^d)$ and $\Pi_m(\R^d)$.
Finally, addressing efficiency, CFs are preferred that only use a small number of data points. 
In particular, one is interested in \emph{interpolatory CFs}, for which $N \leq K$. 
If such a CF also satisfies \ref{item:P1} and \ref{item:P2}, we call it a \emph{positive interpolatory CF}.\footnote{If the cubature weights are nonnegative instead of positive, the CF is called \emph{nonnegative}.}
By now, positive interpolatory---or even minimal---CFs have been constructed for a variety of special cases. 
Most of these are derived for certain (two- and three-dimensional) standard regions, such as the cubes, balls, and triangles. 
Moreover, they usually focus on exactness for certain spaces of algebraic \cite{maxwell1877approximate,stroud1971approximate,mysovskikh1980approximation,cools1993monomial,cools1999monomial,cools2003encyclopaedia} or trigonometric polynomials \cite{cools1997constructing,mysovskikh2001cubature}. 
All of these CFs are derived for certain weight functions and by utilizing specific structures of the integration domain, for instance, being (fully) symmetric.
%
% Tchakaloff 
To the best of my knowledge, the only general result is the following theorem originating from Tchakaloff's work \cite{tchakaloff1957formules} (also see \cite{davis1968nonnegative,bayer2006proof}).
\begin{theorem}[Tchakaloff, 1957]\label{thm:Tchakaloff}
	Given is \eqref{eq:int} with nonnegative $\omega$ and the vector space $\mathcal{F}_K(\Omega)$. 
	Then, there exist $N$ points $\mathbf{y}_1,\dots,\mathbf{y}_N$ in $\Omega$ and positive weights $\lambda_1,\dots,\lambda_N$ with $N \leq K$ such that the corresponding $N$-point CF 
	\begin{equation}\label{eq:ex-cond}
		C_N[f] = \sum_{n=1}^N \lambda_n f(\mathbf{y}_n) 
	\end{equation} 
	is $\mathcal{F}_K(\Omega)$-exact. 
\end{theorem}
In other words: For every compact set $\Omega \subset \R^d$, nonnegative weight function $\omega$, and vector space $\mathcal{F}_K(\Omega)$, there exists a positive interpolatory CF. 
Tchakaloff's original proof is quite beautiful and, in particular, involves the theory of convex bodies. 
Yet, it is not constructive in nature. 
A constructive as well as more elementary proof of Tchakaloff's theorem in the case of $\mathcal{F}_K(\Omega) = \mathbb{P}_m(\R^d)$ and $\omega \equiv 1$ was provided by Davis \cite{davis1968nonnegative}. 
Yet, as it was already noted in \cite{haber1970numerical}, there remain considerable computational difficulties for Davis' approach. 
In fact, these difficulties even increase for other vector spaces $\mathcal{F}_K(\Omega)$ and weight functions $\omega$, and I do not know of any CF constructed from this results.

\subsection*{Novel Contribution}

In the present work, I provide an alternative (constructive) proof of Theorem \ref{thm:Tchakaloff}. 
Yet, the proof presented here is not just more general than the one given in \cite{davis1968nonnegative} but will also be demonstrated to result in the actual construction of positive interpolatory CFs. 
Thereby, only the following restrictions are needed: 
\begin{enumerate}[label=(R\arabic*)] 
	\item \label{item:R3} 
	The integration domain $\Omega \subset \R^d$ is a compact set with positive volume and a boundary of measure zero.

	\item \label{item:R1}
	The vector space $\mathcal{F}_K(\Omega)$ contains constants. In particular, $1 \in \mathcal{F}_K(\Omega)$.
	
	\item \label{item:omega} 
	The weight function $\omega: \Omega \to \R_0^+$ is Riemann integrable\footnote{For a bounded function on a compact set, being Riemann integrable is equivalent to being continuous almost everywhere (in the sense of Lebesgue).} and the set of its zeros ${\{ \, \mathbf{x} \in \Omega \mid \omega(\mathbf{x}) = 0 \, \}}$ is nowhere dense\footnote{A set $A \subset \Omega$ is \emph{nowhere dense} if $\forall \mathbf{x} \in \Omega\setminus A \ \exists \varepsilon > 0 \ \forall \mathbf{a} \in A: \ \|\mathbf{x}-\mathbf{a}\| > \varepsilon$.}. 
	
	\item \label{item:R2} 
	The linear functional $I:\mathcal{F}_K(\Omega) \to \R, f \mapsto I[f]$ is positive definite. 
	
\end{enumerate}

\begin{remark} 
	\ref{item:R3} ensures the existence of sequences that are equidistributed in $\Omega$. 
	\ref{item:R1} and \ref{item:omega} guarantee the construction of nonnegative least squares (LS) CFs.
	In particular, \ref{item:omega} yields the existence of certain discrete inner products. 
	Furthermore, \ref{item:R2} ensures the existence of a continuous inner product. 
\end{remark}

% Approach
The intended procedure to construct positive interpolatory CFs consists of two steps. 
The first step is to construct a nonnegative $N$-point CF that is $\mathcal{F}_K(\Omega)$-exact with $N > K$. 
This is achieved by adapting the LS approach first proposed by Wilson \cite{wilson1970necessary,wilson1970discrete} in one dimension for $\omega \equiv 1$ and $\mathcal{F}_K(\Omega) = \mathbb{P}_{m}(\R)$. 
Later, the approach was further developed by Huybrechs \cite{huybrechs2009stable} and Glaubitz \cite{glaubitz2020shock,glaubitz2020stable,glaubitz2020stableQRs}. 
In a recent work, Glaubitz \cite{glaubitz2020stableCFs} also generalized the LS approach to construct nonnegative exact LS-CFs (for experimental data) in multiple dimensions. 
That is, $d \geq 1$, $\omega \geq 0$, and $\mathcal{F}_K(\Omega) = \mathbb{P}_{m}(\R^d)$. 
Here, I extend these results to arbitrary vector spaces $\mathcal{F}_K(\Omega)$ containing constants. 
Then, the second step consist of utilizing a method of Steinitz \cite{steinitz1913bedingt} (also see \cite{davis1967construction}) to derive a positive interpolatory CF which uses a subset of the LS-CF's data points. 
The Matlab code developed from this procedure can be found at \cite{glaubitz2020github}.

% Advantages and disadvantages
Still, it should be noted that for many standard domains and weight functions there already exist positive interpolatory---or even minimal---CF which are exact on $\mathbb{P}_{m}(\R^d)$ or $\Pi_m(\R^d)$. 
The CFs constructed in the present work are not expected to improve on these formulas. 
Yet, to the best of my knowledge, this work is the first to realize the construction of positive interpolatory CFs for general domains $\Omega$, weight functions $\omega$, and vector spaces $\mathcal{F}_K(\Omega)$. 
These CFs will therefore find their greatest utility away from standard domains and weight functions, or when exactness is not desired for polynomials but for some other vector space of functions (e.\,g.\ radial basis functions or wavelets).
Finally, the present work solves a problem which dates back as far as to Maxwell's original work on approximate multiple integration \cite{maxwell1877approximate} from 1877, and I think its theoretical findings alone might be of high interest for researchers in the field of cubature theory. 

\subsection*{Outline} 

The rest of this work is organized as follows. 
In \S \ref{sec:prelim} some preliminaries on exactness as well as unisolvent and equidistributed sequences are presented. 
Building up on these results, \S \ref{sec:LS} discusses the construction of nonnegative LS-CFs that are exact for a prescribed vector space of functions. 
Next, in \S \ref{sec:procedure}, a procedure to construct positive interpolatory CFs in a quite general setting is proposed. 
Numerical results are presented in \S \ref{sec:num}. 
Finally, some concluding thoughts are offered in \S \ref{sec:summary}.  
\section{Preliminaries: Unisolvent and Equidistributed Sequences} 
\label{sec:prelim} 

Before addressing the generalization of nonnegative LS-CFs (\S \ref{sec:LS}) and the construction of positive interpolatory CFs (\S \ref{sec:procedure}), a few preliminary results are collected. 
These address exactness of CFs as well as unisolvent and equidistributed sequences.

\subsection{Exactness and Unisolvency}

First note that the exactness conditions \eqref{eq:ex-cond} are equivalent to the data points ${X_N = \{ \mathbf{x}_n\}_{n=1}^N}$ and weights of a CF to solve the nonlinear systems \begin{equation}\label{eq:ex-nonlin-system}
  	\underbrace{
  	\begin{pmatrix}
    		\varphi_1(\mathbf{x}_1) & \dots & \varphi_1(\mathbf{x}_N) \\ 
    		\vdots & & \vdots \\ 
    		\varphi_K(\mathbf{x}_1) & \dots & \varphi_K(\mathbf{x}_N)
  	\end{pmatrix}}_{=: \Phi(X_N)}
  	\underbrace{
  	\begin{pmatrix} 
    		w_1 \\ \vdots \\ w_N 
  	\end{pmatrix}}_{=: \mathbf{w}} 
  	= 
  	\underbrace{
  	\begin{pmatrix} 
    		m_1 \\ \vdots \\ m_K 
  	\end{pmatrix}}_{=: \mathbf{m}}.
\end{equation} 
In general, solving \eqref{eq:ex-nonlin-system} can be highly nontrivial, and doing so by brute force may result in a CF which points lie outside of the integration domain or with negative weights \cite{haber1970numerical}. 
The situation changes, however, if a fixed, prescribed, set of data points is assumed. 
Then, $\Phi(X_N) = \Phi$ and \eqref{eq:ex-nonlin-system} becomes a linear system:
\begin{equation}\label{eq:ex-system}
  	\Phi \mathbf{w} = \mathbf{m}
\end{equation} 
The solvability of \ref{eq:ex-system} was studied, for instance, in \cite[Section 2.2]{glaubitz2020stableCFs} for ${\mathcal{F}_K(\Omega) = \mathbb{P}_m(\R^d)}$. 
Here, we extend these results to arbitrary vector spaces $\mathcal{F}_K(\Omega)$ that are spanned by linearly independent continuous functions ${\varphi_1,\dots,\varphi_K:\Omega \to \R}$. 
Note that for $K < N$ \eqref{eq:ex-system} becomes an underdetermined linear system. 
These are well-known to either have no or infinitely many solutions. 
The latter case arises when we restrict ourselves to unisolvent sets of data points. 

\begin{definition}[Unisolvent Point Sets] 
  A set of points $X = \{\mathbf{x}_n\}_{n=1}^N \subset \R^d$ is called \emph{$\mathcal{F}_K(\Omega)$-unisolvent} if for all $f \in \mathcal{F}_K(\Omega)$
  \begin{equation}
    f(\mathbf{x}_n) = 0, \ n=1,\dots,N, \implies f(\boldsymbol{x}) = 0, \ \forall \boldsymbol{x} \in \Omega.
  \end{equation} 
  That is, the only function $f \in \mathcal{F}_K(\Omega)$ that interpolates zero data is the zero function. 
\end{definition} 

Assuming the that $X_N$ is $\mathcal{F}_K(\Omega)$-unisolvent the following result is obtained. 

\begin{lemma}\label{lem:solution-space}
  Let  $K < N$ and $X = \{\mathbf{x}_n\}_{n=1}^N$ be $\mathcal{F}_K(\Omega)$-unisolvent. 
  Then, the linear system \eqref{eq:ex-system} induces an $(N-K)$-dimensional 
affine linear subspace of solutions 
  \begin{equation}\label{eq:sol-space}
    W := \left\{ \mathbf{w} \in \R^N \mid \Phi\mathbf{w}=\mathbf{m} \right\}.
  \end{equation}
\end{lemma} 

\begin{proof} 
	The case $\mathcal{F}_K(\Omega) = \mathbb{P}_m(\R^d)$ was shown in \cite[Lemma 8]{glaubitz2020stableCFs}. 
	It is easy to verify that the same arguments carry over to the general case discussed here. 
\end{proof} 

Next, a simple sufficient criterion for $X_N \subset \Omega$ to be $\mathcal{F}_K(\Omega)$-unisolvent is provided. 
It is based on sequences that are dense in $\Omega \subset \R^d$. 
Recall that $(\mathbf{x}_n)_{n \in \N} \subset \R^d$ is called \emph{dense in $\Omega$} if 
\begin{equation}\label{eq:dense}
	\forall \mathbf{x} \in \Omega \ \, \forall \varepsilon > 0 \ \, \exists n \in \N: \quad 
	\norm{ \mathbf{x} - \mathbf{x}_n } < \varepsilon.
\end{equation}
That is, every point in $\Omega$ can be approximated arbitrarily accurate by an element of the sequence $(\mathbf{x}_n)_{n \in \N}$.\footnote{Note that for the purpose of the present work the definition of density is independent of the norm used in \eqref{eq:dense} since $\Omega$ is located in a finite dimensional space.}  

\begin{lemma}\label{lem:unisolvent}
	Let $(\mathbf{x}_n)_{n \in \N} \subset \Omega$ be dense in $\Omega$ and let $X_N = \{ \mathbf{x}_n \}_{n=1}^N$. 
	Moreover, let $\mathcal{F}_K(\Omega)$ be a vector space spanned by continuous functions $\varphi_1,\dots,\varphi_K:\Omega \to \R$. 
	Then, there exists an $N_0 \in \N$ such that $X_N$ is $\mathcal{F}_K(\Omega)$-unisolvent for every $N \geq N_0$. 
\end{lemma}

\begin{proof}
	Assume that the assertion is wrong. 
	Then, there exists an $f \in \mathcal{F}_K(\Omega)$ with $f \not\equiv 0$ such that ${f(\mathbf{x}_n) = 0}$ for all $n \in \N$. 
	Yet, since $(\mathbf{x}_n)_{n \in \N} \subset \Omega$ is dense in $\Omega$, this either contradicts $f \not\equiv 0$ or $f$ being continuous. 
\end{proof}

If there exists an $N_0 \in \N$ such that $X_N = \{ \mathbf{x}_n \}_{n=1}^N$ is $\mathcal{F}_K(\Omega)$-unisolvent for every $N \geq N_0$, as in Lemma \ref{lem:unisolvent}, $(\mathbf{x}_n)_{n \in \N}$ is referred to as an \emph{$\mathcal{F}_K(\Omega)$-unisolvent} sequence. 
Thus, Lemma \ref{lem:unisolvent} states that every dense sequence is also $\mathcal{F}_K(\Omega)$-unisolvent.

\subsection{Equidistributed Sequences}

Density as a sufficient condition for unisolvency was not discussed in \cite{glaubitz2020stableCFs} but will be handy for the subsequent construction of nonnegative LS-CFs.
Another important property will be for a sequence $(\mathbf{x}_n)_{n \in \N} \subset \Omega$ to satisfy 
\begin{equation}\label{eq:equi} 
	\lim_{N \to \infty} \frac{|\Omega|}{N} \sum_{n=1}^N g(\mathbf{x}_n) 
		= \int_\Omega g(\boldsymbol{x}) \intd \boldsymbol{x} 
\end{equation} 
for all measurable bounded functions $g:\Omega \to \R$ that are continuous almost everywhere (in the sense of Lebesgue). 
It is easy to note that $(\mathbf{x}_n)_{n \in \N}$ simply being dense in $\Omega$ does not suffice to ensure \eqref{eq:equi}.\footnote{It is interesting to note, however, that every dense sequence can be rearranged into a (equidistributed) sequence satisfying \eqref{eq:equi}.}
Yet, in his celebrated work \cite{weyl1916gleichverteilung} from 1916 Weyl showed that \eqref{eq:equi} can be connected to $(\mathbf{x}_n)_{n \in \N}$ being \emph{equidistributed} (also called \emph{uniformly distributed}). 
To not exceed the scope of this work, some details on equidistributed sequences are omit. 
These can be found in the excellent monograph \cite{kuipers2012uniform} of Kuipers and Niederreiter, hoewever.\footnote{In particular, a treatment of compact Hausdorff spaces, as it is the case here, is presented in Chapter 3 of \cite{kuipers2012uniform}.}
Yet, for the purpose of the present work, it is at least worth noting how such an equidistributed sequence $(\mathbf{x}_n)_{n \in \N}$ can be constructed. 

\begin{remark}[Construction of Equidistributed Sequences]\label{eq:eq-points}
	Note that since $\Omega \subset \R^d$ is compact, we can always find an $R > 0$ such that $\Omega$ is contained in the $d$-dimensional hypercube ${C^{(d)}_R = [-R,R]^d}$. 
	For $d=1$, a sequence $(y_n)_{n \in \N}$ that is equidistributed in $[-R,R]$ can be constructed by successively generating a grid of equally spaced points: 
	\begin{equation}
	\begin{aligned}
		y_1 & = -R, \quad y_2 = R, \quad y_3 = 0, \\ 
		y_4 & = -\frac{1}{2}R, \quad y_5 = \frac{1}{2}R, \\ 
		y_6 & = -\frac{3}{4}R, \quad y_7 = -\frac{1}{4}R, \quad 
		y_8 = \frac{1}{4}R, \quad y_9 = \frac{3}{4}R, \quad \dots
	\end{aligned}
	\end{equation} 
	That is, one starts with the end points $y_1 = -R$ and $y_2 = R$. 
	Then, in every step the centers of all already existing pairs of neighboring points are added to the sequence. 
	An equidistributed sequence $(\mathbf{y}_n)_{n \in \N}$ in $C^{(d)}_R$ for $d > 1$ is obtained by a tensor product grid of the above one-dimensional sequence. 
	Finally, a sequence $(\mathbf{x}_n)_{n \in \N}$ that is equidistributed in $\Omega$ is given by the subsequence of $(\mathbf{y}_n)_{n \in \N} \subset C_d(R)$ for which all elements outside of $\Omega$ have been removed, i.\,e.,
	\begin{equation}
		\mathbf{y}_n \in (\mathbf{x}_n)_{n \in \N} \iff \mathbf{y}_n \in \Omega.
	\end{equation} 
	Again, we refer to \cite{kuipers2012uniform} for details. 
	Yet, it is worth noting that for such a constructed sequence $(\mathbf{x}_n)_{n \in \N}$ to be equidistribiuted in $\Omega$ it is necessary for $\Omega$ to have a boundary of measure zero.\footnote{The boundary of $\Omega$ having measure zero ensures that the extension of $g$ to $C^{(d)}_R$ that is zero outside of $\Omega$ is still continuous almost everywhere.} 
	This is ensured by \ref{item:R3}. 
	Another option to construct an euqidistributed sequence, which was also used for the three-dimensional numerical test in \S \ref{sec:num}, are Halton points \cite{halton1960efficiency} (a generalization of the one-dimensional van der Corput points \cite[Erste Mitteilung]{van1935verteilungsfunktionen}).\footnote{Both sequences belong to the family of low-discrepancy points, which are developed to minimize the upper bound provided by the famous Koksma--Hlawak inequality \cite{hlawka1961funktionen,niederreiter1992random}.}
\end{remark} 

Finally, it should be stressed that equidistributed sequences are dense sequences with a specific ordering. 
This preliminary section is therefore closed by the following corollary. 

\begin{corollary}\label{cor:equid}
	Let $\Omega \subset \R^d$ be a compact set and let its boundary have measure zero. 
	Furthermore, let $(\mathbf{x}_n)_{n \in \N}$ be the sequence described in Remark \ref{eq:eq-points}. 
	Then, $(\mathbf{x}_n)_{n \in \N}$ is equidistributed in $\Omega$. 
	In particular, it satisfies \eqref{eq:equi} and is $\mathcal{F}_K(\Omega)$-unisolvent. 
\end{corollary}  
\section{Least Squares Cubature Formulas}
\label{sec:LS} 

In this section, it is demonstrated how nonnegative $\mathcal{F}_K(\Omega)$-exact LS-CFs can be constructed by using a sufficiently large set of data points coming from an equidistributed sequence. 
This is done by generalizing the LS approach from previous works \cite{wilson1970necessary,wilson1970discrete,huybrechs2009stable,glaubitz2020shock,glaubitz2020stableQRs,glaubitz2020stableCFs}. 
In particular, the results of \cite{glaubitz2020stableCFs} are extended from $\mathbb{P}_m(\Omega)$ to general vector spaces $\mathcal{F}_K(\Omega)$.

\subsection{Formulation as an LS Problem}
\label{sub:LS-problem}

Let $(\mathbf{x}_n)_{n \in \N}$ be a $\mathcal{F}_K(\Omega)$-unisolvent sequence in $\Omega$ and let $X_N = \{\mathbf{x}_n\}_{n=1}^N$. 
As noted before, an $\mathcal{F}_K(\Omega)$-exact CF with data points $X_N$ can be constructed by determining a vector of cubature weights that solves the linear system of exactness conditions \eqref{eq:ex-system}. 
For $N > K$, \eqref{eq:ex-system} induces an $(N-K)$-dimensional affine linear subspace space of solutions $W \subset \R^N$. 
All of these yield an $\mathcal{F}_K(\Omega)$-exact $N$-point CF.
The LS approach consists of finding the unique vector of weights $\mathbf{w} \in W$ that minimizes a weighted Euclidean norm: 
\begin{equation}\label{eq:LS-sol} 
	\mathbf{w}^{\mathrm{LS}} = \argmin_{\mathbf{w} \in W} \ \norm{ R^{-1/2} \mathbf{w} }_{2}
\end{equation} 
This vector is called the \emph{LS solution} of $\Phi \mathbf{w} = \mathbf{m}$. 
Here, $R^{-1/2}$ is a diagonal \emph{weight matrix}, given by 
\begin{equation} 
	R^{-1/2} = \diag\left( \frac{1}{\sqrt{r_1}}, \dots, \frac{1}{\sqrt{r_N}} \right), \quad 
	r_n > 0, \quad n=1,\dots,N.
\end{equation} 
If $r_n = 0$, this is interpreted as a constraint $w_n = 0$.
The corresponding $N$-point CF 
\begin{equation}\label{eq:LS-CF}
	C^{\mathrm{LS}}_N[f] = \sum_{n=1}^N w_n^{\mathrm{LS}} f(\mathbf{x}_n)
\end{equation}
is called an \emph{LS-CF}. 
In \S \ref{sub:positivity} it is shown that choosing the \emph{discrete weights} $r_n$ as 
\begin{equation} 
	r_n = \frac{\omega(\mathbf{x}_n) |\Omega|}{N}, \quad n=1,\dots,N,
\end{equation} 
results in the LS-CFs to be not only $\mathcal{F}_K(\Omega)$-exact but also conditionally nonnegative ($N$ needs to be sufficiently large).
% Explicit representation
Finally, it is convenient to note that---at least formally---the LS solution is explicitly given by 
\begin{equation}\label{eq:LS-sol2}
  \mathbf{w}^{\mathrm{LS}} = R \Phi^T (\Phi R \Phi^T)^{-1} \mathbf{m};   
\end{equation} 
see \cite{cline1976l_2}.
Thereby, $R P^T (P R P^T)^{-1}$ is the Moore--Penrose pseudoinverse of $R^{-1/2}P$; see \cite{ben2003generalized}.
By utilizing a beautiful connection to discrete orthonormal bases, \eqref{eq:LS-sol2} can be considerably simplified; see \S \ref{sub:char}.

\subsection{Orthonormal Bases}
\label{sub:ON-bases}

Recall that the linear functional $I:\mathcal{F}_K(\Omega) \to \R, f \mapsto I[f]$ is assumed to be positive definite; see \ref{item:R2}.
Hence, 
\begin{equation}\label{eq:cont-inner-prod}
  \scp{u}{v} = \int_{\Omega} u(\boldsymbol{x}) v(\boldsymbol{x}) \omega(\boldsymbol{x}) \intd \boldsymbol{x}, \quad 
  \norm{u} = \sqrt{\scp{u}{u}}
\end{equation} 
define an inner product and a corresponding norm on $\mathcal{F}_K(\Omega)$.
In particular, \eqref{eq:cont-inner-prod} induces an orthonormal basis $\{ \pi_k \}_{k=1}^K$ of $\mathcal{F}_K(\Omega)$ then. 
That is, the functions $\pi_1,\dots,\pi_K$ span the space $\mathcal{F}_K(\Omega)$ while satisfying 
\begin{equation}
  \scp{\pi_k}{\pi_l} = \delta_{k,l} := 
  	\begin{cases} 
		1 &; \ k = l, \\ 
		0 &; \ k \neq l,
	\end{cases} 
	\quad k,l = 1,\dots,K.
\end{equation} 
Such a basis is referred to as a \emph{continuous orthonormal basis (COB)} and its elements are denoted by $\pi_k(\cdot;\omega)$.
%
% dicrete inner product 
Analogously, assuming that ${X_N^+ = \{ \, \mathbf{x}_n \mid r_n > 0, \ n=1,\dots,N \, \}}$ is $\mathcal{F}_K(\Omega)$-unisolvent,
\begin{equation}\label{eq:disc-inner-prod}
	[u,v]_N = \sum_{n=1}^N r_n u(\mathbf{x}_n) v(\mathbf{x}_n), \quad 
	\norm{u}_N = \sqrt{[u,u]_N}
\end{equation} 
defines a discrete inner product and a corresponding norm on $\mathcal{F}_K(\Omega)$. 
Of course, also \eqref{eq:disc-inner-prod} induces an orthonormal basis. 
In contrast to the COB, this basis satisfies 
\begin{equation}
	[\pi_k,\pi_l]_N = \delta_{k,l}, \quad k,l = 1,\dots,K,
\end{equation} 
while spanning $\mathcal{F}_K(\Omega)$. 
It is therefore referred to as a \emph{discrete orthonormal basis (DOB)} and its elements are denoted by $\pi_k(\cdot;\mathbf{r})$. 
% 
% Gram--Schmidt
Henceforth, it is assumed that both bases, the COB and DOB, are constructed by (modified) Gram--Schmidt orthonormalization applied to the same initial basis $\{\varphi_k\}_{k=1}^K$:
\begin{equation}\label{eq:Gram-Schmidt}
\begin{aligned}
  \tilde{\pi}_k(\boldsymbol{x};\omega) & = \varphi_k(\boldsymbol{x}) - \sum_{l=1}^{k-1} \scp{\varphi_k}{\pi_l(\cdot;\omega)} \pi_l(\boldsymbol{x};\omega), \quad 
  && \pi_k(\boldsymbol{x};\omega) = \frac{\tilde{\pi}_k(\boldsymbol{x};\omega)}{\norm{\tilde{\pi}_k(\cdot;\omega)}}, \\ 
  \tilde{\pi}_k(\boldsymbol{x};\mathbf{r}) & = \varphi_k(\boldsymbol{x}) - \sum_{l=1}^{k-1} [\varphi_k,\pi_l(\cdot;\mathbf{r})]_N \pi_l(\boldsymbol{x};\mathbf{r}), \quad 
  && \pi_k(\boldsymbol{x};\mathbf{r}) = \frac{\tilde{\pi}_k(\boldsymbol{x};\mathbf{r})}{\norm{\tilde{\pi}_k(\cdot;\mathbf{r})}_N},
\end{aligned}
\end{equation} 
for $\boldsymbol{x} \in \Omega$. 
See one of the excellent textbooks \cite{gautschi1997numerical,trefethen1997numerical,gautschi2004orthogonal} or \cite{glaubitz2020shock,glaubitz2020stableCFs}.
Furthermore, it shall be assumed that the vector space $\mathcal{F}_K(\Omega)$ contains constants; see \ref{item:R1}.
Hence, w.\,l.\,o.\,g., $\varphi_1 \equiv 1$ and therefore 
\begin{equation} 
	\pi_1(\cdot;\omega) \equiv \frac{1}{\norm{1}}, \quad 
	\pi_1(\cdot;\mathbf{r}) \equiv \frac{1}{\norm{1}_N}. 
\end{equation}
This will turn out to be convenient to prove conditionally nonnegativity of the LS-CFs in \S \ref{sub:positivity}.

\subsection{Characterization of the LS Solution}
\label{sub:char}

In my opinion, the real beauty of the LS approach is revealed when combined with the concept of DOBs. 
Observe that the matrix product $\Phi R \Phi^T$ in the explicit representation of the LS solution \eqref{eq:LS-sol2} can be interpreted as a Gram matrix w.\,r.\,t.\ the discrete inner product \eqref{eq:disc-inner-prod}: 
\begin{equation}
  \Phi R \Phi^T = 
  \begin{pmatrix}
    [\varphi_1,\varphi_1]_N & \dots & 
    [\varphi_1,\varphi_K]_N \\ 
    \vdots & & \vdots \\ 
    [\varphi_K,\varphi_1]_N & \dots & 
    [\varphi_K,\varphi_K]_N \\ 
  \end{pmatrix}
\end{equation}
Thus, if the linear system \eqref{eq:ex-system} is formulated w.\,r.\,t.\ $\{\varphi_k(\cdot ; \mathbf{r})\}_{k=1}^K$, one gets $\Phi R \Phi^T = I$. 
This further yields \eqref{eq:LS-sol2} to become  
\begin{equation}\label{eq:LS-sol3}
	\mathbf{w}^{\mathrm{LS}} = R \Phi^T \mathbf{m}.  
\end{equation}
In particular, the LS weights are explicitly given by 
\begin{equation}\label{eq:LS-sol-explicit}
	w_n^{\mathrm{LS}} = r_n \sum_{k=1}^K \pi_k( \mathbf{x}_n ; \mathbf{r}) I[ \pi_k( \, \cdot \, ; \mathbf{r} ) ], \quad n=1,\dots,N. 
\end{equation}
Subsequently, this representation allows me to prove that the LS-CFs are conditionally nonnegative.

\subsection{Nonnegativity of LS-CFs}
\label{sub:positivity} 

Fist, two technical lemmas are given. 
Both lemmas were proven in \cite{glaubitz2020stableCFs} for the case $\mathcal{F}_K(\Omega) = \mathbb{P}_m(\Omega)$. 
It is easy to verify that they also hold for a general vector space of continuous functions and their proofs are therefore omit. 

\begin{lemma}\label{lem:1}
	Let $\Omega \subset \R^d$ be a compact set and assume that 
	\begin{equation}\label{eq:assum-lemma}
		\lim_{N \to \infty} [u,v]_N = \scp{u}{v} \quad \forall u,v \in \mathcal{F}_K(\Omega).
	\end{equation} 
	Furthermore, let $(u_N)_{N \in \N}, (u_N)_{N \in \N} \subset \mathcal{F}_K(\Omega)$ and $u, v \in \mathcal{F}_K(\Omega)$ such that 
	\begin{equation} 
		\lim_{N \to \infty} u_N = u, \quad 
		\lim_{N \to \infty} u_N = u \quad 
		\text{in } ( \mathcal{F}_K(\Omega), \norm{\cdot}_\infty ).
	\end{equation} 
	Then, 
	\begin{equation} 
		\lim_{N \to \infty} [u_N,v_N]_N = \scp{u}{v}.
	\end{equation}
\end{lemma}

Here, $\norm{\cdot}_\infty$ denotes the usual maximum norm with $\norm{f}_\infty = \max_{\boldsymbol{x} \in \Omega} |f(\boldsymbol{x})|$.

\begin{lemma}\label{lem:2}
	Let $\Omega \subset \R^d$ be a compact set, let $\{ \varphi_k \}_{k=1}^K$ be a basis of $\mathcal{F}_K(\Omega)$, and assume that \eqref{eq:assum-lemma} holds. 
	Moreover, let $\{ \pi_k(\cdot,\omega) \}_{k=1}^K$ and $\{ \pi_k(\cdot,\mathbf{r}) \}_{k=1}^K$ respectively denote the COB and DOB constructed from $\{ \varphi_k \}_{k=1}^K$ by Gram-Schmidt orthonormalization \eqref{eq:Gram-Schmidt}. 
	Then, 
	\begin{equation} 
		\lim_{N \to \infty} \pi_k(\cdot,\mathbf{r}) = \pi_k(\cdot,\omega) \quad 
		\text{in } ( \mathcal{F}_K(\Omega), \norm{\cdot}_\infty ).
	\end{equation} 
\end{lemma} 

Note that $\mathbf{r} = (r_1,\dots,r_N)$ and the DOB therefore depends on $N \in \N$. 
Lemmas \ref{lem:1} and \ref{lem:2} essentially state that if the discrete inner product converges to the continuous inner product, then also the DOB converges (uniformly) to the COB. 
Finally, this observation enables one to proof the following theorem. 

\begin{theorem}[The LS-CF is conditionally nonnegative]\label{thm:main}
	Let $\Omega \subset \R^d$ be a compact set and $1\in \mathcal{F}_K(\Omega)$. 
	Moreover, let $(\mathbf{x}_n)_{n \in \N} \subset \Omega$ and $(r_n)_{n \in \N} \subset \R_0^+$ such that  
 	\begin{equation}\label{eq:cond-thm}
		\lim_{N \to \infty} \sum_{n=1}^N r_n u(\mathbf{x}_n) v(\mathbf{x}_n) 
			= \int_\Omega \omega(\boldsymbol{x}) u(\boldsymbol{x}) v(\boldsymbol{x}) \intd \boldsymbol{x} 
			\quad \forall u, v \in \mathcal{F}_K(\Omega). 
	\end{equation} 
	Furthermore, assume that the subsequence $(\mathbf{x}^+_n)_{n \in \N} \subset (\mathbf{x}_n)_{n \in \N}$ containing only the $\mathbf{x}_n$ with $r_n > 0$ is $\mathcal{F}_K(\Omega)$-unisolvent.
	Then, there exists an $N_0 \in \N$ such that for all $N \geq N_0$ the corresponding LS-CF \eqref{eq:LS-CF} is nonnegative. 
\end{theorem} 

\begin{proof} 
	First, it should be noted that $(\mathbf{x}^+_n)_{n \in \N}$ being $\mathcal{F}_K(\Omega)$-unisolvent ensures the existence of a discrete inner product and therefore of a DOB. 
	W.\,l.\,o.\,g., it can be assumed that this DOB, $\{ \pi_k(\cdot,\mathbf{r}) \}_{k=1}^K$, as well as the COB,  $\{ \pi_k(\cdot,\omega) \}_{k=1}^K$, are constructed by applying Gram-Schmidt orthonormalization to the same initial basis $\{ \varphi_k \}_{k=1}^K$. 
	Hence, the LS weights are explicitly given by \eqref{eq:LS-sol-explicit} for $N \geq K$. 
	Next, let 
	\begin{equation} 
		\epsilon_k := [\pi_k(\cdot,\mathbf{r}),1]_N - \scp{\pi_k(\cdot,\mathbf{r})}{1}. 
	\end{equation} 
	Then, the LS weights can be rewritten as 
	\begin{equation}
		w_n^{\mathrm{LS}}  
			= r_n \left( \pi_1(\mathbf{x}_n;\mathbf{r}) [\pi_1(\cdot;\mathbf{r}),1]_N 
			- \sum_{k=1}^K \varepsilon_k \pi_k(\mathbf{x}_n;\mathbf{r}) \right).
	\end{equation} 
	Note that, w.\,l.\,o.\,g., $\varphi \equiv 1$ and $\pi_1(\cdot,\mathbf{r}) \equiv 1/\|1\|_N$. 
	This yields 
	\begin{equation}
    		\pi_1(\mathbf{x}_n;\mathbf{r}) [\pi_1(\cdot;\mathbf{r}),1]_N = 1. 
  	\end{equation} 
	If $r_n = 0$, it directly follows that $w_n^{\text{LS}} = 0$.
	For $r_n > 0$, on the other hand, the assertion $w_n^{\text{LS}} \geq 0$ is equivalent to 
	\begin{equation}\label{eq:assertion2-thm}
    		\sum_{k=1}^K \varepsilon_k \pi_k(\mathbf{x}_n;\mathbf{r}) \leq 1.
  	\end{equation} 
	At the same time, \eqref{eq:cond-thm} and Lemma \ref{lem:2} imply that every element of the DOB converges uniformly to the corresponding element of the COB. 
	In particular, for every $k=1,\dots,K$, the sequence ${(\pi_k(\cdot,\mathbf{r}))_{N \in \N} \subset \mathcal{F}_K(\Omega)}$ is uniformly bounded. 
	Thus, there exists a constant $C > 0$ such that 
	\begin{equation}
    		\sum_{k=1}^K \varepsilon_k \pi_k(\mathbf{x}_n;\mathbf{r}) 
      		\leq C \sum_{k=1}^K \left| \varepsilon_k \right|.
  	\end{equation} 
	Moreover, Lemma \ref{lem:1} implies $\lim_{N \to \infty} \epsilon_k \to 0$ for all  $k=1,\dots,K$. 
	Hence, there exists an $N_0 \geq K$ such that 
	\begin{equation}
    		\left| \varepsilon_k \right| \leq \frac{1}{CK}, \quad k=1,\dots,K,
  	\end{equation}
  	for $N \geq N_0$. 
  	Finally, this yields \eqref{eq:assertion2-thm} and therefore the assertion. 
\end{proof}

A simple consequence of Theorem \ref{thm:main} is the subsequent corollary in which the special case of the equidistributed points from Remark \ref{eq:eq-points} is considered. 

\begin{corollary}\label{cor:main}
	Given are $\Omega \subset \R^d$, $\omega: \Omega \to \R_0^+$, and $\mathcal{F}_K(\Omega) \subset C(\Omega)$ such that the restrictions \ref{item:R3}--\ref{item:R2} are satisfied. 
	Let $(\mathbf{x}_n)_{n \in \N} \subset \Omega$ be the equidistributed sequence as in Remark \ref{eq:eq-points}. 
	Then, there exists an $N_0 \in \N$ such that for all $N \geq N_0$ and discrete weights 
	\begin{equation} 
		r_n = \frac{|\Omega| \omega(\mathbf{x}_n)}{N}, \quad n=1,\dots,N,
	\end{equation} 
	the corresponding LS-CF \eqref{eq:LS-CF} is $\mathcal{F}_K(\Omega)$-exact as well as nonnegative. 
\end{corollary} 

\begin{proof} 
	It has been shown in Corollary \ref{cor:equid} that, under the restrictions \ref{item:R3} and \ref{item:omega}, the sequence $(\mathbf{x}_n)_{n \in \N} \subset \Omega$ is $\mathcal{F}_K(\Omega)$-unisolvent and satisfies \eqref{eq:equi} for all measurable bounded functions that are continuous almost everywhere. 
	In particular, \eqref{eq:cond-thm} holds for $(\mathbf{x}_n)_{n \in \N}$ and the discrete weights $(r_n)_{n \in \N}$. 
	Furthermore, \ref{item:omega} ensures the subsequence $(\mathbf{x}^+_n)_{n \in \N}$ to be $\mathcal{F}_K(\Omega)$-unisolvent as well. 
	In combination with \ref{item:R1} and \ref{item:R2}, Theorem \ref{thm:main} implies the assertion. 
\end{proof} 

Building up on this result, subsequently, a simple procedure to construct positive interpolatory CFs is developed.   
\section{Proposed Procedure} 
\label{sec:procedure} 

In this section, the previous theoretical findings from \S \ref{sec:prelim} and \S \ref{sec:LS} are forged into a rigorous procedure to construct positive interpolatory CFs. 
This is done under the assumption that $\Omega \subset \R^d$, $\omega: \Omega \to \R_0^+$, and $\mathcal{F}_K(\Omega) \subset C(\Omega)$ satisfy the restrictions \ref{item:R3}--\ref{item:R2}.
The proposed procedure consists of two steps: 

\begin{enumerate}[label=(S\arabic*)] 
	\item \label{item:S1}
	Construct an $\mathcal{F}_K(\Omega)$-exact nonnegative $N$-point CF utilizing the LS approach. 
	
	\item \label{item:S2} 
	If $N>K$, use Steinitz' method to successively reduce the number of data points until $N \leq K$. 
	
\end{enumerate}

The result will be a positive interpolatory CF that uses a subset of at most $K$ data points. 
The realization of this procedure is demonstrated in \S \ref{sec:num}.

\subsection{Constructing Positive LS-CFs} 

Given is the equidistributed sequence ${(\mathbf{x}_n)_{n \in \N} \subset \Omega}$ as in Remark \ref{eq:eq-points}. 
The basic idea is to increase the set of data points ${X_N = \{\mathbf{x}_n\}_{n=1}^N}$ until the corresponding LS-CF is nonnegative. 
One starts with $N = K$ and proceeds as follows:  
(1) If $X_N$ is $\mathcal{F}_K(\Omega)$-unisolvent, the discrete weights are chosen as $r_n = |\Omega|\omega(\mathbf{x}_n)/N$, and (2) the corresponding LS weights $\mathbf{w}^{\mathrm{LS}}$ are computed. 
If these weights are all nonnegative, a nonnegative LS-CF has been constructed. 
Otherwise, meaning that $X$ is not $\mathcal{F}_K(\Omega)$-unisolvent or $\mathbf{w}^{\mathrm{LS}}$ is not nonnegative, the number of data points $N$ is, for instance, doubled and one returns to (1). 
This procedure is summarized below in Algorithm \ref{algo:LS-CF}. 

\begin{algorithm}
\caption{Construction of a nonnegative LS-CF}
\label{algo:LS-CF}
\begin{algorithmic}[1]
    \State{$N = K$, $r = 0$, and $w_{\text{min}} = 0$} 
    \While{$r < K$ or $w_{\text{min}} < 0$} 
      	\State{$X_N = \{\mathbf{x}_n\}_{n=1}^N$} 
		\State{Compute a DOB, for instance, by the modified Gram-Schmidt procedure}
      	\State{Compute the matrix $\Phi = \Phi(X_N)$}
      	\State{Compute the rank of $\Phi$: $r = \text{rank}(\Phi)$} 
      	\If{$r = K$}
      		\State{Compute the LS weights $\mathbf{w}^{\text{LS}}$ as in \eqref{eq:LS-sol3}} 
      		\State{Determine the smallest weight: $w_{\text{min}} = \min( \mathbf{w}^{\text{LS}} )$}
		\EndIf 
	\State{$N = 2N$} 
    \EndWhile 
\end{algorithmic}
\end{algorithm} 

Note that $K = r := \text{rank}(\Phi)$ is equivalent to $X_N$ being $\mathcal{F}_K(\Omega)$-unisolvent. 
It should also be stressed that,  by Corollary \ref{cor:main}, Algorithm \ref{algo:LS-CF} is ensured to terminate.  
The output is an $\mathcal{F}_K(\Omega)$-exact nonnegative LS-CF using $N$ data points. 
Usually, we can expect $N$ to be larger than $K$. 
In what follows, I therefore adapt a method of Steinitz' \cite{steinitz1913bedingt} (also see \cite{davis1967construction}) to successively reduce the number of data points from $N$ to $K$

\subsection{Reducing the Number of Data Points: Steinitz' Method} 
\label{sub:Steinitz} 

Given is an $\mathcal{F}_K(\Omega)$-exact nonnegative LS-CF. 
In a first step, all cubature weights that are equal to zero as well as the corresponding data points are removed. 
This yields an $\mathcal{F}_K(\Omega)$-exact \emph{positive} CF 
\begin{equation}\label{eq:pos-CF}
	C_N[f] = \sum_{n=1}^N w_n f(\mathbf{x}_n),
\end{equation}
where $N$ denotes the number of data points in a generic sense.  
If still $N> K$, one proceeds by using a method of Steinitz \cite{steinitz1913bedingt} (also see \cite{davis1967construction}). 
Here, this method is referred to as \emph{Steinitz' method}.
It allows one to successively reduce the number of data points until $N \leq K$ by going over to an appropriate subset of data points, while preserving $\mathcal{F}_K(\Omega)$-exactness as well as nonnegativity. 

Recall that the vector space $\mathcal{F}_K(\Omega)$ has dimension $K$, and so does its algebraic dual space (the space of all linear functionals defined on $\mathcal{F}_K(\Omega)$). 
In particular, among the $N$ linear functionals 
\begin{equation} 
	L_n[f] = f(\mathbf{x}_n), \quad n=1,\dots,N,
\end{equation}
at most $K$ are linearly independent. 
That is, if $N > K$, there exist a vector of coefficients $\mathbf{a} = (a_1, \dots, a_M)^T$ such that 
\begin{equation}\label{eq:vec-a}
	a_1 L_1[f] + \dots + a_N L_N[f] = 0 \quad \forall f \in \mathcal{F}_K(\Omega)
\end{equation}
and $a_n > 0$ for at least one $n$.
Let $\sigma$ be given by  
\begin{equation}
	\sigma = \max_{1 \leq n \leq N} \frac{a_n}{w_n}. 
\end{equation}
Then, $\sigma > 0$, $\sigma w_n - a_n \geq 0$ for all $n$, and $\sigma w_n - a_n = 0$ for at least one $n$. 
From \eqref{eq:vec-a} one therefore has 
\begin{equation}
	I[f] = \frac{\sigma w_1 - a_1}{\sigma} L_1[f] + \dots + \frac{\sigma v_N - a_N}{\sigma} L_N[f] 
	\quad \forall f \in \mathcal{F}_K(\Omega).
\end{equation} 
Note that one of the coefficients is zero and, together with the corresponding linear functional (data point), can be removed.
Hence, on $\mathcal{F}_K$, $I$ has been expressed as a linear combination of at most $N-1$ of the linear functionals $L_1,\dots,L_N$ with positive coefficients. 

Iterating this process of removing all zero weights and then applying Steinitz' method, one finally arrives at a positive interpolatory CF 
\begin{equation}
	C_N[f] = \sum_{k=1}^N \lambda_k f(\mathbf{y}_k) 
\end{equation} 
with $N \leq K$. 
That is, the CF $C_N$ is $\mathcal{F}_K(\Omega)$-exact while using at most $K$ data points $\mathbf{y}_k$---all of which lie inside of $\Omega$---and only positive cubature weights $\lambda_k$. 
The whole procedure is summarized in Algorithm \ref{algo:Steinitz}

\begin{algorithm}
\caption{Reduce the number of data points}
\label{algo:Steinitz}
\begin{algorithmic}[1]
	\While{$K < N$} 
    		\State{Compute $\Phi = \Phi(X)$ and $\text{null}(\Phi)$} 
      	\State{Determine $\mathbf{a} \in \text{null}(\Phi) \setminus \{\mathbf{0}\}$ s.\,t.\ $a_n > 0$ for at least one $n$} 
		\State{Compute $\sigma = \max_{n} a_n/w_n$} 
		\State{Overwrite the cubature weights: $w_n = (\sigma w_n - a_n)/\sigma$} 
		\State{Remove all zero weights as well as the corresponding data points} 
		\State{$N = N -  \#\{ \, w_n \mid w_n = 0, \ n=1,\dots,N \, \}$} 
    \EndWhile 
\end{algorithmic}
\end{algorithm} 

Here, $\text{null}(\Phi)$ denotes the null space of the matrix $\Phi$. 
For $\Phi \in \R^{K \times N}$, it is defined as 
\begin{equation}
	\text{null}(\Phi) = \{ \, \mathbf{a} \in \R^N \mid \Phi \mathbf{a} = \mathbf{0} \, \}
\end{equation}
and yields a linear subspace of $\R^N$.

\begin{remark}[Construction of $\mathbf{a}$]
	Note that \eqref{eq:vec-a} is equivalent to $\mathbf{a} \in \text{null}(\Phi)$. 
	Thus, essentially every $\mathbf{a} \in \text{null}(\Phi) \setminus \{\mathbf{0}\}$ will do the job (if $a_n \leq 0$ for all $n=1,\dots,M$, one can go over to $-a$). 
	Moreover, it was shown in Lemma \ref{lem:solution-space} that $\text{null}(\Phi)$ has dimension $N-K$. 
	Hence, as long as $N > K$, such a vector of coefficients $\mathbf{a}$ can always be found. 
\end{remark}
\section{Numerical Results} 
\label{sec:num}

In this section some numerical result for the proposed positive interpolatory CFs are presented. 
The corresponding Matlab code can be found at \cite{glaubitz2020github}.
Focus is given to the cases ${\mathcal{F}_K(\Omega) = \mathbb{P}_m(\R^d)}$ and ${\mathcal{F}_K(\Omega) = \Pi_m(\R^d)}$. 
For the algebraic polynomials, the \emph{total degree} ${|\boldsymbol{\alpha}| = \sum_{j=1}^d} |\alpha_j|$ for ${\boldsymbol{\alpha} = (\alpha_1,\dots,\alpha_d)}$ is considered. 
Other choices for the degree are possible and include the \emph{absolute degree} and the \emph{Euclidean degree}; see \cite{haber1970numerical,cools1997constructing,trefethen2017cubature}. 
Trigonometric polynomials, on the other hand, are linear combinations $\sum_{|\boldsymbol{\alpha}| \leq m} c_{\boldsymbol{\alpha}} t_{\boldsymbol{\alpha}}$ of trigonometric monomials 
\begin{equation}
	\prod_{j=1}^d e^{2 \pi i \alpha_j x_j} 
	\quad \text{with} \quad i^2 = -1,
\end{equation} 
where $\boldsymbol{\alpha} \in \mathbb{Z}$. 
Furthermore, it is required that $c_{\boldsymbol{\alpha}}$ and $c_{-\boldsymbol{\alpha}}$ are complex conjugates so that all the monomials are real-valued. 
In this case, one has 
\begin{equation}
	\text{dim} \ \mathbb{P}_m(\R^d) = \text{dim} \ \Pi_m(\R^d) = \binom{m+d}{m}
\end{equation} 
for the dimension of both vector spaces. 
Yet, it should be stressed once more that the specific choice of the vector space does not affect the success of the proposed procedure.

%%% Tables %%%

\begin{table}[!htb]
\renewcommand{\arraystretch}{1.3}
\centering 
%\begin{adjustbox}{width=1\textwidth}
  	\begin{tabular}{c c c c c c c c c c c c c c c c} 
	\toprule
	\multicolumn{14}{c}{Test Case A} \\ \hline
	\multicolumn{4}{c}{$m=0$ ($K=1$)} & & \multicolumn{4}{c}{$m=1$ ($K=4$)} & & \multicolumn{4}{c}{$m=2$ ($K=6$)} \\ \hline 
	$n$ & $x_n$ & $y_n$ & $w_n$ & $\quad$ & $n$ & $x_n$ & $y_n$ & $w_n$ & $\quad$ & $n$ & $x_n$ & $y_n$ & $w_n$ \\ \hline 
	 1 & 0 & 0 & 4 &  & 1 & 1/3 & -1/3 & 4/3 &   & 1 & -1    & 0.2  & 0.8 \\ 
	    &    &    &    &  & 2 & 1/3 & 1     & 4/3 &   & 2 & -0.6 & -1    & 0.8 \\ 
	    &    &    &    &  & 3 & 1    & 1/3  & 4/3 &   & 3 & -0.6 & -0.6 & 0.8 \\ 
	    &    &    &    &  &   &  &  &  &                   & 4 & -0.2 & -0.2 & 0.494... \\ 
	    &    &    &    &  &   &  &  &  &                   & 5 & -0.2 & 0.2  & 0.305... \\ 
	    &    &    &    &  &   &  &  &  &                   & 6 & 0.2  & 0.6  & 0.8 \\ \bottomrule
\end{tabular} 
%\end{adjustbox}
\caption{Data points and cubature weights for $\Omega = C^{(2)}$, $\omega \equiv 1$, and $\mathcal{F}_K(\Omega) = \Pi_m(\R^2)$}
\label{tab:C2-trig}
\end{table}

\begin{table}[!htb]
\renewcommand{\arraystretch}{1.3}
\centering 
\begin{adjustbox}{width=1\textwidth}
  	\begin{tabular}{c c c c c c c c c c c c c c c c c c c} 
	\toprule
	\multicolumn{17}{c}{Test Case B} \\ \hline
	\multicolumn{5}{c}{$m=0$ ($K=1$)} & & \multicolumn{5}{c}{$m=1$ ($K=4$)} & & \multicolumn{5}{c}{$m=2$ ($K=10$)} \\ \hline 
	$n$ & $x_n$ & $y_n$ & $z_n$ & $w_n$ & $\quad$ & $n$ & $x_n$ & $y_n$ & $z_n$ & $w_n$ & $\quad$ & $n$ & $x_n$ & $y_n$ & $z_n$ & $w_n$ \\ \hline 
	 1 & 0 & 0 & -1 & $I[1]$ &  & 1 & 0 & 0  & -1 & 0.598... &  & 1   & -0.5 & 0     & -0.5 & 0.815... \\ 
	 	&    &   &      &           &  & 2 & 0 & -1 & 0  & 1.196... &  & 2   & 0.5  & -0.5 & -0.5 & 0.190... \\ 
	 	&    &   &      &           &  & 3 & 0 & 1  & 0  & 1.196... &  & 3   & 0.5  & 0.5  & -0.5 & 0.516... \\ 
		&    &   &      &           &  & 4 & 0 & 0  & 1  & 0.598... &  & 4   & 0     & -1    & 0     & 0.326... \\ 
	 	&    &   &      &           &  &    &    &     &     &              &  & 5   & 0.5  & 0     & 0     & 0.217... \\ 
	 	&    &   &      &           &  &    &    &     &     &              &  & 6   & -0.5 & 0     & 0.5  & 0.489... \\ 
	 	&    &   &      &           &  &    &    &     &     &              &  & 7   & -0.5 & 0.5  & 0.5  & 0.217... \\ 
	 	&    &   &      &           &  &    &    &     &     &              &  & 8   & 0     & 0.5  & 0.5  & 0.217... \\ 
	 	&    &   &      &           &  &    &    &     &     &              &  & 9   & 0.5  & -0.5 & 0.5  & 0.353... \\ 
	 	&    &   &      &           &  &    &    &     &     &              &  & 10 & 0.5  & 0.5  & 0.5  & 0.244... \\  \bottomrule
\end{tabular} 
\end{adjustbox}
\caption{Data points and cubature weights for $\Omega = B^{(3)}$, $\omega(\boldsymbol{x}) = \sqrt{ \|\boldsymbol{x}\|_2}$, and $\mathcal{F}_K(\Omega) = \mathbb{P}_m(\R^3)$}
\label{tab:B3-alg}
\end{table}

\begin{table}[!htb]
\renewcommand{\arraystretch}{1.3}
\centering 
%\begin{adjustbox}{width=1\textwidth}
  	\begin{tabular}{c c c c c c c c c c c c c c c c} 
	\toprule
	\multicolumn{14}{c}{Test Case C} \\ \hline
	\multicolumn{4}{c}{$m=0$ ($K=1$)} & & \multicolumn{4}{c}{$m=1$ ($K=3$)} & & \multicolumn{4}{c}{$m=2$ ($K=6$)} \\ \hline 
	$n$ & $x_n$ & $y_n$ & $w_n$ & $\quad$ & $n$ & $x_n$ & $y_n$ & $w_n$ & $\quad$ & $n$ & $x_n$ & $y_n$ & $w_n$ \\ \hline 
	 1 & 0 & 0 & I[1] &  & 1 & -2/3 & 2/3 & 1.695... &  & 1 & -2/3 & 2/3  & 0.562... \\ 
	    &    &    &       &  & 2 & 2/3 & -2/3 & 1.695... &  & 2 & 0     & -2/3 & 1.115... \\ 
	    &    &    &       &  & 3 & 2    & 2     & 0.75      &  & 3 & 0     & 0     & 0.928... \\ 
	    &    &    &       &  &    &       &        &              &  & 4 & 0     & 2/3  & 0.267... \\ 
	    &    &    &       &  &    &       &        &              &  & 5 & 4/3  & 4/3  & 0.991... \\ 
	    &    &    &       &  &    &       &        &              &  & 6 & 2     & 4/3  & 0.276... \\ \bottomrule
\end{tabular} 
%\end{adjustbox}
\caption{Data points and cubature weights for $\Omega = B^{(2)} \cup C^{(2)}_{1/2}(1.5,1.5)$, $\omega \equiv 1$, and $\mathcal{F}_K(\Omega) = \mathbb{P}_m(\R^2)$}
\label{tab:combi-alg}
\end{table}

Tables \ref{tab:C2-trig}, \ref{tab:B3-alg} and \ref{tab:combi-alg} list the data points and cubature weights of some specific positive interpolatory CFs constructed by the procedure proposed in \S \ref{sec:procedure}.
Table \ref{tab:C2-trig} addresses the case where $\omega \equiv 1$ and $\Omega$ is the two-dimensional hypercube $C^{(2)}$ that is centered at $(0,0)$ and with radius $1$. 
Generally speaking, the following definition are used:  
\begin{equation} 
\begin{aligned}
	C^{(d)}_r(\mathbf{x}_0) 
		& := \{ \, \mathbf{x} \in \R^d \mid \norm{\mathbf{x}-\mathbf{x}_0}_{\infty} \leq r \, \}, \quad 
	C^{(d)} = C^{(d)}_1(\mathbf{0}), \\ 
	B^{(d)}_r(\mathbf{x}_0) 
		& := \{ \, \mathbf{x} \in \R^d \mid \norm{\mathbf{x}-\mathbf{x}_0}_2 \leq r \, \}, \quad 
	B^{(d)} = B^{(d)}_1(\mathbf{0}) 
\end{aligned}
\end{equation}  
Moreover, $\mathcal{F}_K(\Omega) = \Pi_m(\R^2)$ is considered. 
It is demonstrated in Table \ref{tab:C2-trig} that the number of data points is always less or equal to the number of linearly independent basis functions which are treated exactly by the CF.  
At the same time, all data points lie inside of $\Omega$ and the cubature weights are all positive. 
Table \ref{tab:B3-alg} reports the same for the three-dimensional ball $B^{(3)}$ with nonconstant weight function $\omega(\boldsymbol{x}) = \sqrt{ \|\boldsymbol{x}\|_2}$. 
In this case, however, $\mathcal{F}_K(\Omega) = \mathbb{P}_m(\R^3)$ is considered.
Finally, Table \ref{tab:combi-alg} addresses a nonstandard domain $\Omega \subset \R^2$ given as the union of $B^{(2)}$ and $C_{1/2}^{(2)}(1.5,1.5)$. 
Furthermore, $\omega \equiv 1$ and $\mathcal{F}_K(\Omega) = \mathbb{P}_m(\R^2)$. 

%%% Illustration for non-standard domain %%%

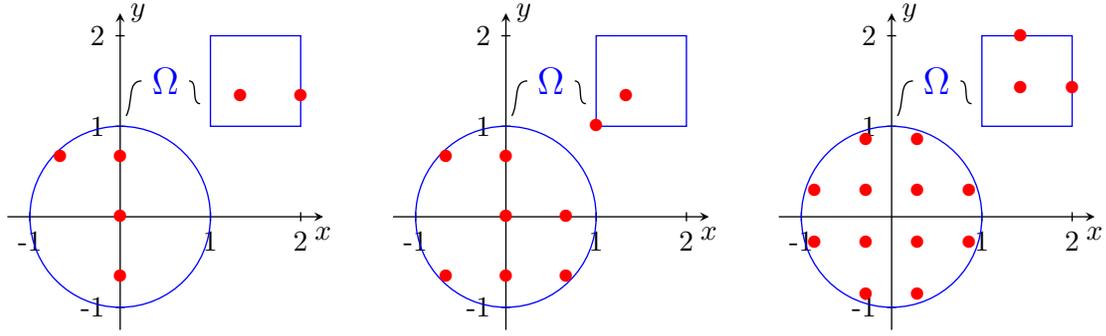
\begin{figure}[!htb]
	\centering 
	\begin{subfigure}[b]{0.32\textwidth} 
    \begin{center}
  	\begin{tikzpicture}[domain = -2.5:4.5, scale=.6, line width=0.5pt]

			% x- and y-axis 
			\draw[->,>=stealth] (-2.5,0) -- (4.5,0) node[below] {$x$};
    			\draw[->,>=stealth] (0,-2.5) -- (0,4.5) node[right] {$y$};
			% labels 
			\draw (-2,0.1) -- (-2,-0.1) node [below] {-1};
			\draw (2,0.1) -- (2,-0.1) node [below] {1};
			\draw (4,0.1) -- (4,-0.1) node [below] {2};
			\draw (-0.1,-2) -- (0.1,-2) node [left] {-1 \ };
			\draw (-0.1,2) -- (0.1,2) node [left] {1 \ };
			\draw (-0.1,4) -- (0.1,4) node [left] {2 \ };

			% Circle 
			\draw[blue]  (0,0) circle [radius = 2]; 
			% Rectangle 
			\draw[blue] (2,2) rectangle (4,4);
			% add label 
			\draw (0,2) node (A) {}; % node on circle
			\draw (2,2.5) node (B) {}; % node on rectangle
			\draw[blue] (1,3) node (C) {\Large $\Omega$};
			\draw (A) to[out=60, in=180] (C);
			\draw (B) to[out=180, in=0] (C);

			% draw data points 
			\foreach \Point in { (-1.33,1.33), (0,-1.33), (0,0), (0,1.33), (2.66,2.66), (4,2.66) }{
    				\node[red] at \Point {\large{\textbullet}};
			}

	\end{tikzpicture}
	\end{center}
	\caption{$m = 2$ ($K=6$) and $N=6$}
    \label{fig:points-TestC-m2} 
  	\end{subfigure}
	\begin{subfigure}[b]{0.32\textwidth} 
    \begin{center}
  	\begin{tikzpicture}[domain = -2.5:4.5, scale=.6, line width=0.5pt]

			% x- and y-axis 
			\draw[->,>=stealth] (-2.5,0) -- (4.5,0) node[below] {$x$};
    			\draw[->,>=stealth] (0,-2.5) -- (0,4.5) node[right] {$y$};
			% labels 
			\draw (-2,0.1) -- (-2,-0.1) node [below] {-1};
			\draw (2,0.1) -- (2,-0.1) node [below] {1};
			\draw (4,0.1) -- (4,-0.1) node [below] {2};
			\draw (-0.1,-2) -- (0.1,-2) node [left] {-1 \ };
			\draw (-0.1,2) -- (0.1,2) node [left] {1 \ };
			\draw (-0.1,4) -- (0.1,4) node [left] {2 \ };

			% Circle 
			\draw[blue]  (0,0) circle [radius = 2]; 
			% Rectangle 
			\draw[blue] (2,2) rectangle (4,4);
			% add label 
			\draw (0,2) node (A) {}; % node on circle
			\draw (2,2.5) node (B) {}; % node on rectangle
			\draw[blue] (1,3) node (C) {\Large $\Omega$};
			\draw (A) to[out=60, in=180] (C);
			\draw (B) to[out=180, in=0] (C);

			% draw data points 
			\foreach \Point in { (-1.33,-1.33), (-1.33,1.33), (0,-1.33), (0,0), (0,1.33), (1.33,-1.33) , (1.33,0), (2.66,2.66), (2,2)}{
    				\node[red] at \Point {\large{\textbullet}};
			}

	\end{tikzpicture}
	\end{center}
	\caption{$m = 3$ ($K=10$) and $N=9$}
    \label{fig:points-TestC-m3} 
  	\end{subfigure}
	\begin{subfigure}[b]{0.32\textwidth} 
    \begin{center}
  	\begin{tikzpicture}[domain = -2.5:4.5, scale=0.6, line width=0.5pt]

			% x- and y-axis 
			\draw[->,>=stealth] (-2.5,0) -- (4.5,0) node[below] {$x$};
    			\draw[->,>=stealth] (0,-2.5) -- (0,4.5) node[right] {$y$};
			% labels 
			\draw (-2,0.1) -- (-2,-0.1) node [below] {-1};
			\draw (2,0.1) -- (2,-0.1) node [below] {1};
			\draw (4,0.1) -- (4,-0.1) node [below] {2};
			\draw (-0.1,-2) -- (0.1,-2) node [left] {-1 \ };
			\draw (-0.1,2) -- (0.1,2) node [left] {1 \ };
			\draw (-0.1,4) -- (0.1,4) node [left] {2 \ };

			% Circle 
			\draw[blue]  (0,0) circle [radius = 2]; 
			% Rectangle 
			\draw[blue] (2,2) rectangle (4,4);
			% add label 
			\draw (0,2) node (A) {}; % node on circle
			\draw (2,2.5) node (B) {}; % node on rectangle
			\draw[blue] (1,3) node (C) {\Large $\Omega$};
			\draw (A) to[out=60, in=180] (C);
			\draw (B) to[out=180, in=0] (C);

			% draw data points 
			\foreach \Point in { (-1.71,-0.57), (-1.71,0.57), (-0.57,-1.71), (-0.57,-0.57), (-0.57,0.57), (-0.57,1.71), (0.57,-1.71), (0.57,-0.57), (0.57,0.57), (0.57,1.71), (1.71,-0.57), (1.71,0.57), (2.85,2.85), (2.85,4), (4,2.85)}{
    				\node[red] at \Point {\large{\textbullet}};
			}

	\end{tikzpicture}
	\end{center}
	\caption{$m = 4$ ($K=15$) and $N=15$}
    \label{fig:points-TestC-m4} 
  	\end{subfigure}
  	\caption{Illustration of the data points for $\Omega = B^{(2)} \cup C^{(2)}_{1/2}(1.5,1.5)$, $\omega \equiv 1$, and $\mathcal{F}_K(\Omega) = \mathbb{P}_m(\R^2)$}
  	\label{fig:points}
\end{figure}

This nonstandard domain $\Omega$, together with the data points of the corresponding positive interpolatory CF, are further illustrated in Figure \ref{fig:points} for $m=2,3,4$. 
In particular, it can be noted that in some cases even $N < K$.

%%% Accuracy %%%

\begin{figure}[!htb]
\centering
	\begin{subfigure}[b]{0.32\textwidth}
		\includegraphics[width=\textwidth]{%
      		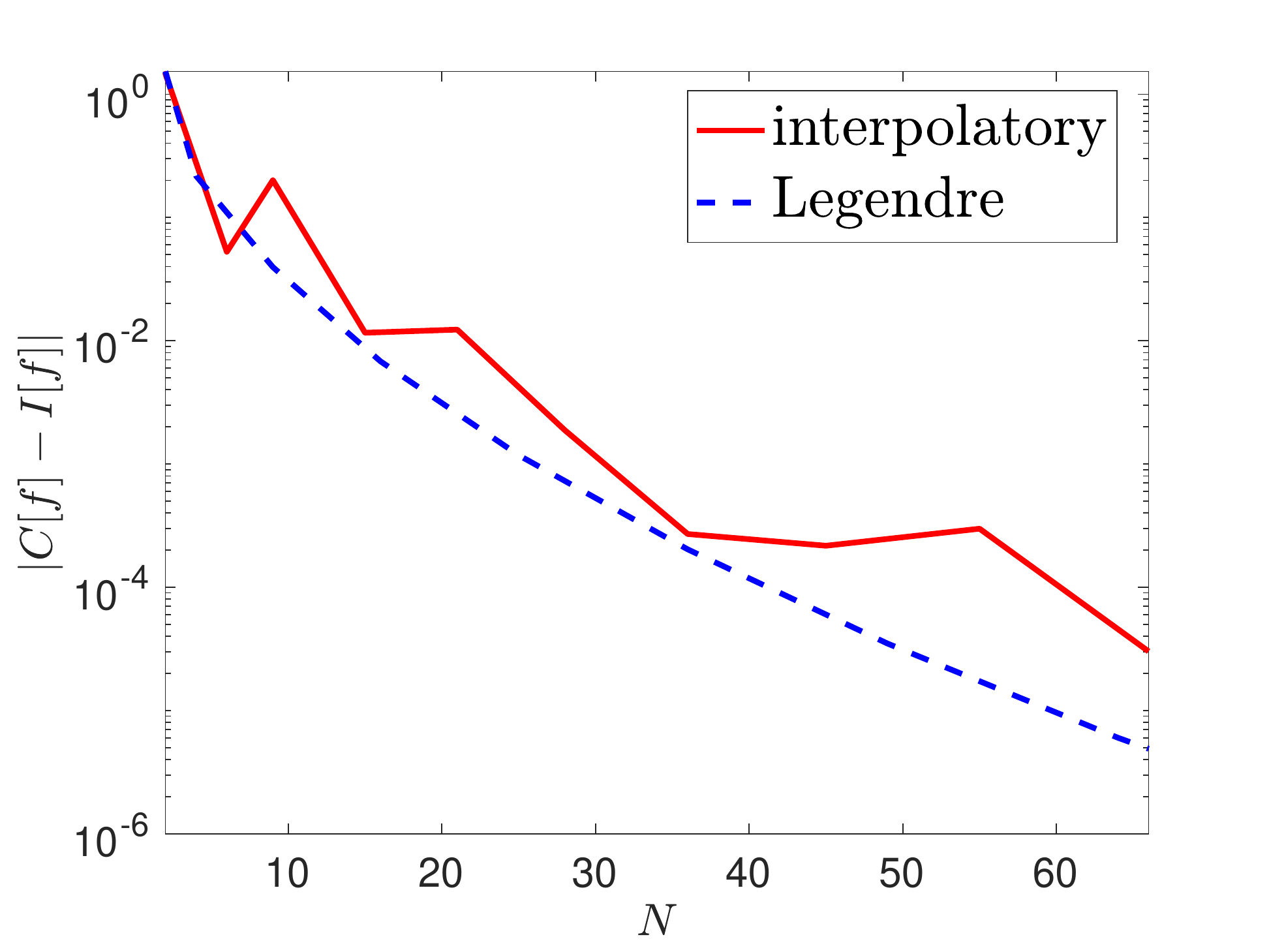}  
    		\caption{$C^{(2)}$ with $\omega \equiv 1$}
    		\label{fig:accuracy-C2-1}
	\end{subfigure}%
	~
	\begin{subfigure}[b]{0.32\textwidth}
    		\includegraphics[width=\textwidth]{%
      		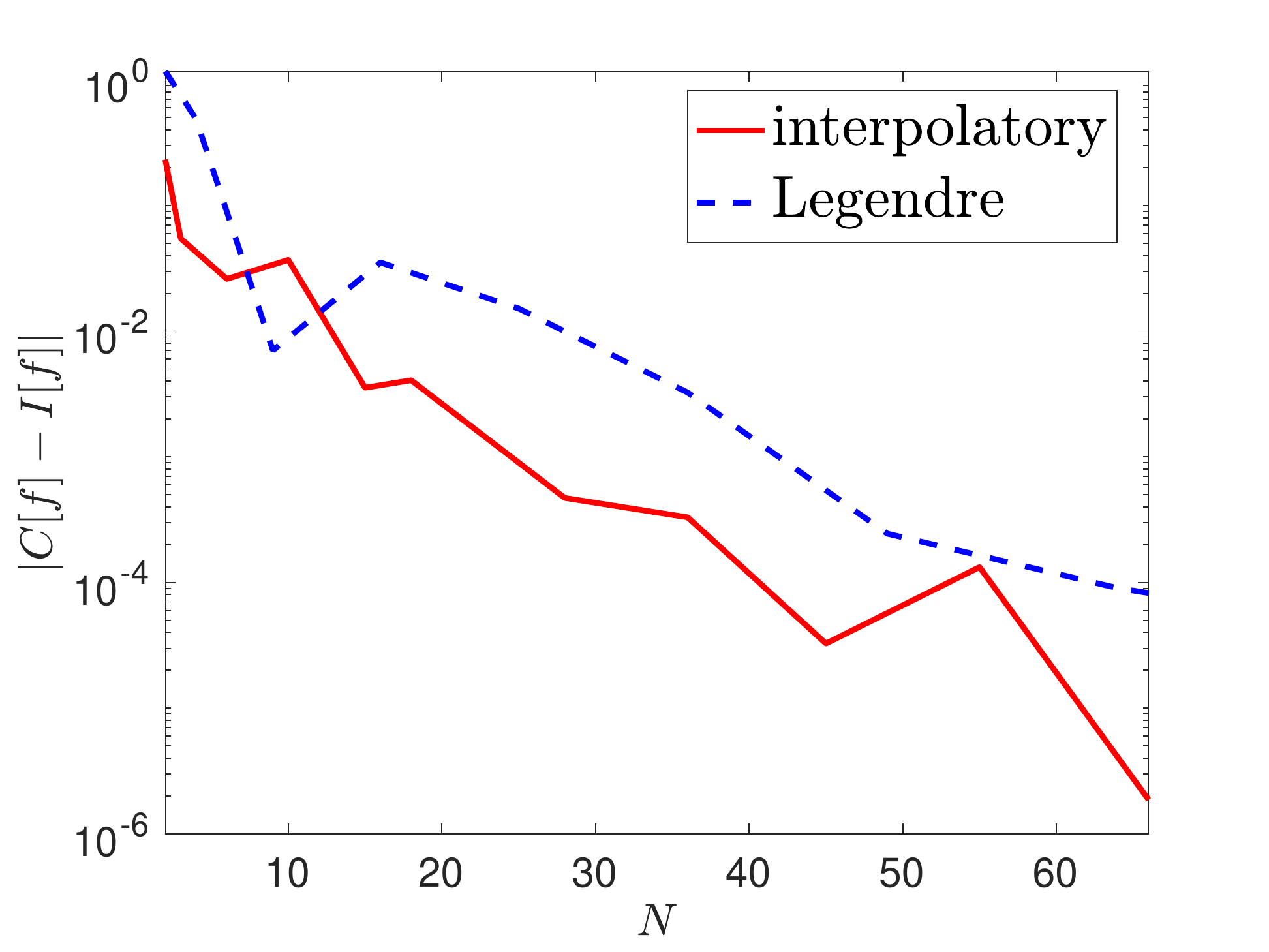} 
    		\caption{$B^{(2)}$ with $\omega \equiv 1$}
    		\label{fig:accuracy-B2-1}
  	\end{subfigure}%
  	~
  	\begin{subfigure}[b]{0.32\textwidth}
   		\includegraphics[width=\textwidth]{%
      		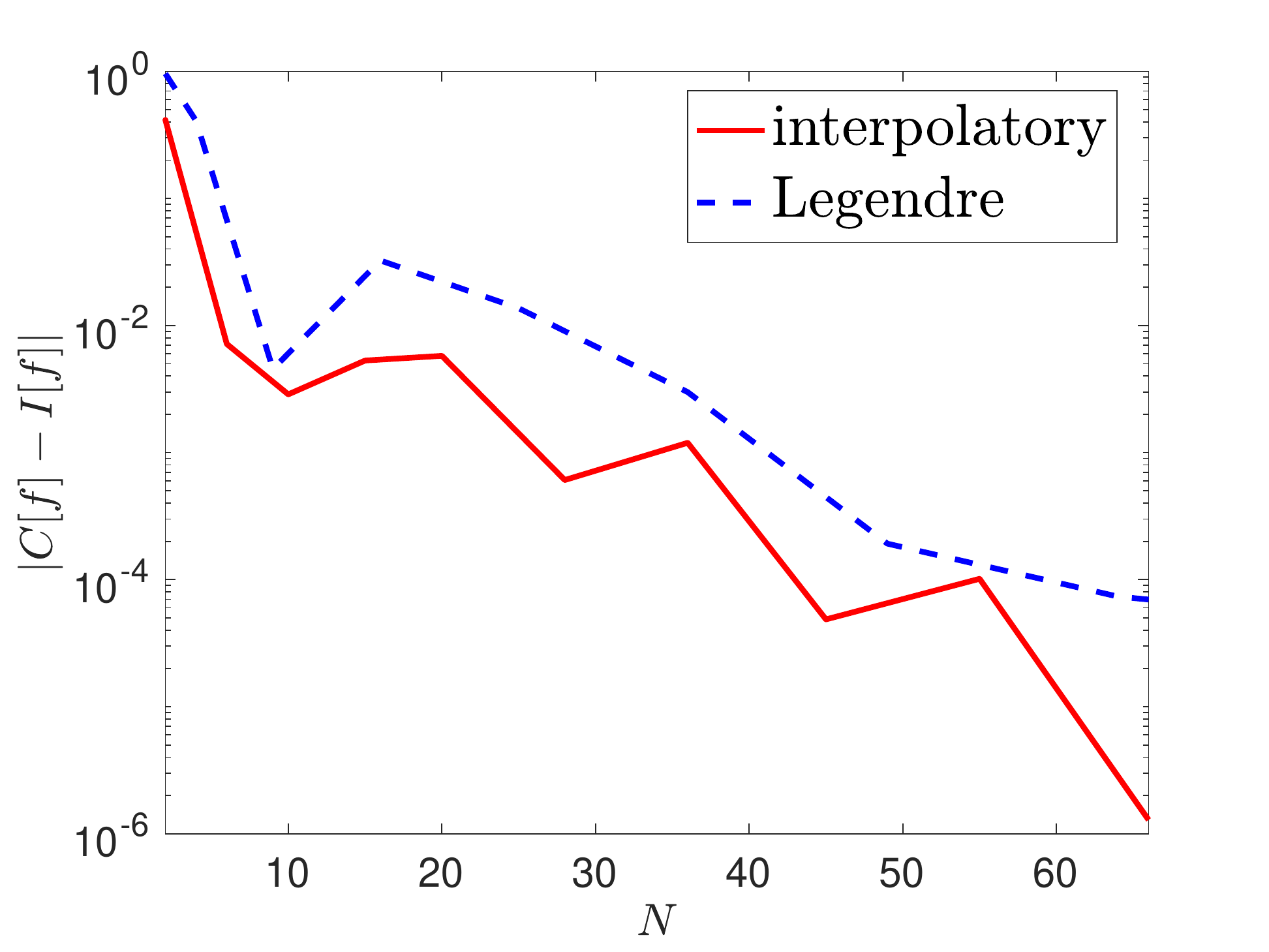} 
    		\caption{$B^{(2)}$ with $\omega(\boldsymbol{x}) = \sqrt{\|\boldsymbol{x}\|_2}$}
    		\label{fig:accuracy-B2-sqrt}
  	\end{subfigure}%
  	\\ 
  	\begin{subfigure}[b]{0.32\textwidth}
    		\includegraphics[width=\textwidth]{%
      		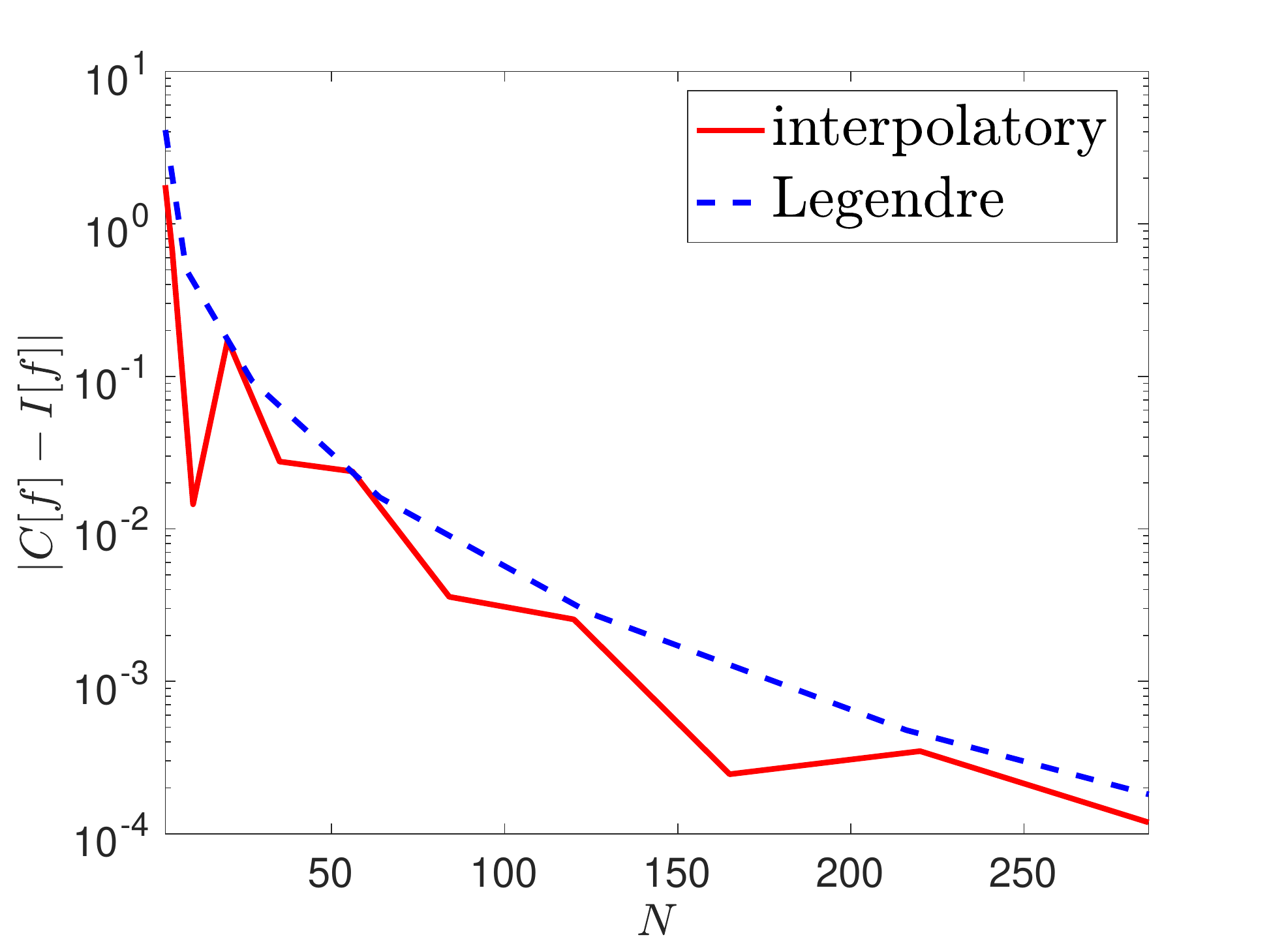} 
    		\caption{$C^{(3)}$ with $\omega \equiv 1$}
    		\label{fig:accuracy-C3-1}
  	\end{subfigure}%
  	~
  	\begin{subfigure}[b]{0.32\textwidth}
    		\includegraphics[width=\textwidth]{%
      		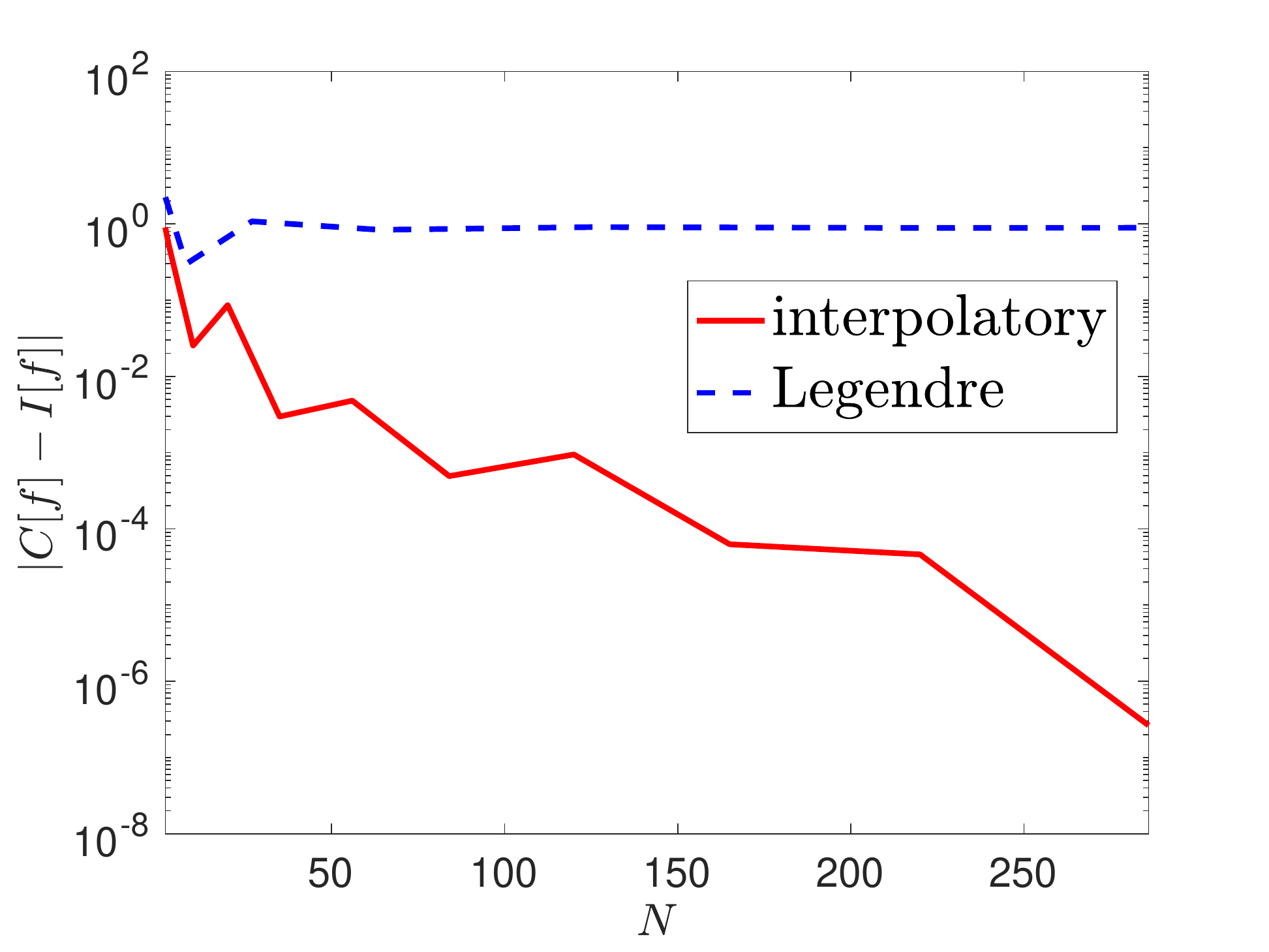} 
    		\caption{$B^{(3)}$ with $\omega \equiv 1$}
    		\label{fig:accuracy-B3-1}
  	\end{subfigure}%
  	~
  	\begin{subfigure}[b]{0.32\textwidth}
    		\includegraphics[width=\textwidth]{%
      		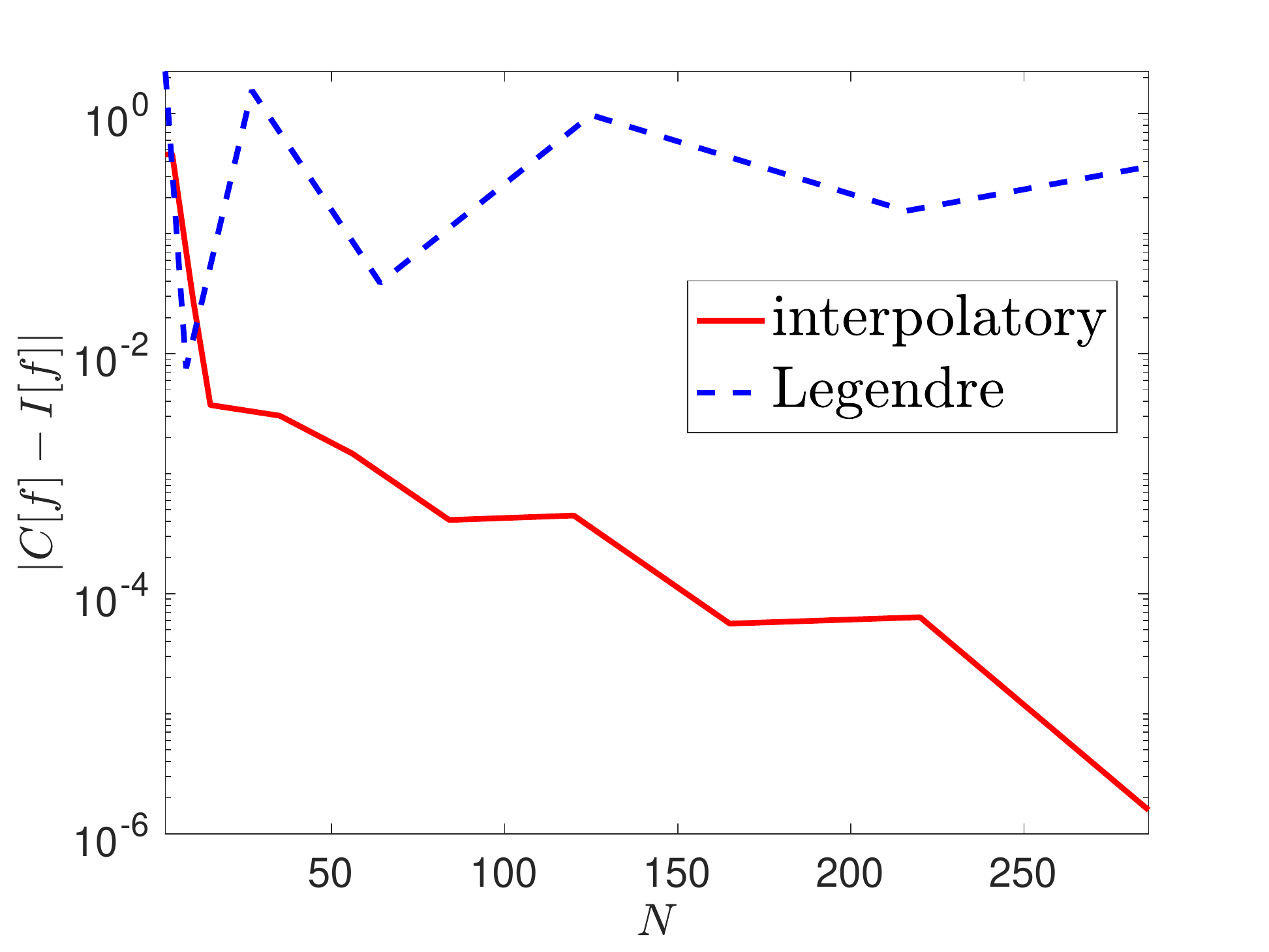} 
    		\caption{$B^{(3)}$ with $\omega(\boldsymbol{x}) = \sqrt{\|\boldsymbol{x}\|_2}$}
    		\label{fig:accuracy-B3-sqrt}
  	\end{subfigure}%
  	\caption{Accuracy of the constructed positive interpolatory CF compared to Gauss--Legendre-type rules. 
	In figures \ref{fig:accuracy-C2-1} and \ref{fig:accuracy-C3-1}, $f(\boldsymbol{x}) = f(x_1)\dots f(x_d)$ with $f(x) = 1/(1+x^2)$ is considered. 
	In all other figures, $f(\boldsymbol{x}) = 1/(1+\|\boldsymbol{x}\|_2^2) + \sin(x_1)$ is used. 
	The positivie interpolatory CF was constructed to be exact for $\mathcal{F}_K(\Omega) = \mathbb{P}_m(\R^d)$ with increasing $m$.}
	\label{fig:accuracy}
\end{figure}

Finally, a comparison of the positive interpolatory CFs with Gauss--Legendre-type rules is provided. 
Figure \ref{fig:accuracy} reports on the accuracy of both rules for six different test cases. 
For $\Omega = B^{(d)}$ with $d=2,3$, a transformed product Gauss--Legendre rule has been used; see \cite{davis2007methods}. 
It should be stressed that the transformation of this rule to $B{(d)}$ is not exact for algebraic polynomials anymore. 
As a result, especially for $\Omega = B{(3)}$, the constructed positive interpolatory CFs have a significantly higher accuracy. 
Yet, also in all other cases, the positive interpolatory CFs are able to keep up with the Gauss--Legendre-type rules. 
In nonstandard regions, such as unions of different (standard) domains, the positive interpolatory CFs can be expected to perform even better compared to known CFs---assuming there even exists a known CF. 
\section{Summary} 
\label{sec:summary} 

In this work, a procedure was developed to construct positive interpolatory CFs in a fairly general setting. 
In particular, only the restrictions \ref{item:R3}--\ref{item:R2} are needed. 
While the existence of such CFs has theoretically been proven already in 1957 by Tchakaloff, their actual construction was not achieved until now. 
The present work fills this gap in the theory of cubature. 
% future work
%In a future work, I aim to develop a second---more elegant and even more general---method to construct positive interpolatory CFs by leveraging a connection between positive interpolatory CFs and sparse signal recovery. 

%\section*{Acknowledgements}
%This work has benefited from helpful advice from ...

\bibliographystyle{siamplain}
\bibliography{literature}

\begin{thebibliography}{10}

\bibitem{bayer2006proof}
{\sc C.~Bayer and J.~Teichmann}, {\em The proof of {T}chakaloff’s theorem},
  Proceedings of the AMS, 134 (2006), pp.~3035--3040.

\bibitem{ben2003generalized}
{\sc A.~Ben-Israel and T.~N. Greville}, {\em Generalized Inverses: Theory and
  Applications}, vol.~15 of CMS Books in Mathematics, Springer Science \&
  Business Media, 2003.

\bibitem{boyer2011history}
{\sc C.~B. Boyer and U.~C. Merzbach}, {\em A History of Mathematics}, John
  Wiley \& Sons, 2011.

\bibitem{cline1976l_2}
{\sc R.~Cline and R.~J. Plemmons}, {\em $\ell_2$-solutions to underdetermined
  linear systems}, SIAM Review, 18 (1976), pp.~92--106.

\bibitem{cools1997constructing}
{\sc R.~Cools}, {\em Constructing cubature formulae: The science behind the
  art}, Acta Numerica, 6 (1997), pp.~1--54.

\bibitem{cools1999monomial}
{\sc R.~Cools}, {\em Monomial cubature rules since “{S}troud”: a
  compilation—part 2}, Journal of Computational and Applied Mathematics, 112
  (1999), pp.~21--27.

\bibitem{cools2003encyclopaedia}
{\sc R.~Cools}, {\em An encyclopaedia of cubature formulas}, Journal of
  Complexity, 19 (2003), pp.~445--453.

\bibitem{cools1993monomial}
{\sc R.~Cools and P.~Rabinowitz}, {\em Monomial cubature rules since
  “{S}troud”: a compilation}, Journal of Computational and Applied
  Mathematics, 48 (1993), pp.~309--326.

\bibitem{davis1968nonnegative}
{\sc P.~Davis and M.~Wilson}, {\em Nonnegative interpolation formulas for
  uniformly elliptic equations}, Journal of Approximation Theory, 1 (1968),
  pp.~374--380.

\bibitem{davis1967construction}
{\sc P.~J. Davis}, {\em A construction of nonnegative approximate quadratures},
  Mathematics of Computation, 21 (1967), pp.~578--582.

\bibitem{davis2007methods}
{\sc P.~J. Davis and P.~Rabinowitz}, {\em Methods of Numerical Integration},
  Courier Corporation, 2007.

\bibitem{engels1980numerical}
{\sc H.~Engels}, {\em Numerical Quadrature and Cubature}, Academic Press, 1980.

\bibitem{gautschi1997numerical}
{\sc W.~Gautschi}, {\em Numerical Analysis}, Springer Science \& Business
  Media, 1997.

\bibitem{gautschi2004orthogonal}
{\sc W.~Gautschi}, {\em Orthogonal Polynomials: Computation and Approximation},
  Oxford University Press, 2004.

\bibitem{glaubitz2020github}
{\sc J.~Glaubitz}, {\em jglaubitz/positive\_interpolatory\_{CF}s}, 2020,
  \url{https://doi.org/10.5281/zenodo.4019333}.

\bibitem{glaubitz2020shock}
{\sc J.~Glaubitz}, {\em Shock Capturing and High-Order Methods for Hyperbolic
  Conservation Laws}, Logos Verlag Berlin GmbH, 2020.

\bibitem{glaubitz2020stableCFs}
{\sc J.~Glaubitz}, {\em Stable high-order cubature formulas for experimental
  data},  (2020).
\newblock Submitted.

\bibitem{glaubitz2020stableQRs}
{\sc J.~Glaubitz}, {\em Stable high order quadrature rules for scattered data
  and general weight functions}, SIAM Journal on Numerical Analysis, 58 (2020),
  pp.~2144--2164.

\bibitem{glaubitz2020stable}
{\sc J.~Glaubitz and P.~{\"O}ffner}, {\em Stable discretisations of high-order
  discontinuous {G}alerkin methods on equidistant and scattered points},
  Applied Numerical Mathematics, 151 (2020), pp.~98--118.

\bibitem{haber1970numerical}
{\sc S.~Haber}, {\em Numerical evaluation of multiple integrals}, SIAM Review,
  12 (1970), pp.~481--526.

\bibitem{halton1960efficiency}
{\sc J.~H. Halton}, {\em On the efficiency of certain quasi-random sequences of
  points in evaluating multi-dimensional integrals}, Numerische Mathematik, 2
  (1960), pp.~84--90.

\bibitem{hlawka1961funktionen}
{\sc E.~Hlawka}, {\em Funktionen von beschr{\"a}nkter {V}ariation in der
  {T}heorie der {G}leichverteilung}, Ann. Mat. Pura Appl., 54 (1961),
  pp.~325--333.

\bibitem{huybrechs2009stable}
{\sc D.~Huybrechs}, {\em Stable high-order quadrature rules with equidistant
  points}, Journal of Computational and Applied Mathematics, 231 (2009),
  pp.~933--947.

\bibitem{kuipers2012uniform}
{\sc L.~Kuipers and H.~Niederreiter}, {\em Uniform Distribution of Sequences},
  Courier Corporation, 2012.

\bibitem{maxwell1877approximate}
{\sc J.~C. Maxwell}, {\em On approximate multiple integration between limits of
  summation}, in Proc. Cambridge Philos. Soc, vol.~3, 1877, pp.~39--47.

\bibitem{mysovskikh1980approximation}
{\sc I.~Mysovskikh}, {\em The approximation of multiple integrals by using
  interpolatory cubature formulae}, in Quantitative Approximation, Elsevier,
  1980, pp.~217--243.

\bibitem{mysovskikh2001cubature}
{\sc I.~P. Mysovskikh}, {\em Cubature formulae that are exact for trigonometric
  polynomials}, TW Reports,  (2001).
\newblock Edited by R. Cools and H.J. Schmid.

\bibitem{niederreiter1992random}
{\sc H.~Niederreiter}, {\em Random Number Generation and Quasi-Monte Carlo
  Methods}, SIAM, 1992.

\bibitem{steinitz1913bedingt}
{\sc E.~Steinitz}, {\em Bedingt konvergente {R}eihen und konvexe {S}ysteme},
  Journal f{\"u}r die reine und angewandte Mathematik (Crelle’s Journal),
  1913 (1913), pp.~128--176.

\bibitem{stroud1971approximate}
{\sc A.~H. Stroud}, {\em Approximate Calculation of Multiple Integrals},
  Prentice-Hall, 1971.

\bibitem{tchakaloff1957formules}
{\sc V.~Tchakaloff}, {\em Formules de cubatures m{\'e}caniques {\`a}
  coefficients non n{\'e}gatifs}, Bull. Sci. Math, 81 (1957), pp.~123--134.

\bibitem{trefethen2017cubature}
{\sc L.~N. Trefethen}, {\em Cubature, approximation, and isotropy in the
  hypercube}, SIAM Review, 59 (2017), pp.~469--491.

\bibitem{trefethen1997numerical}
{\sc L.~N. Trefethen and D.~Bau~III}, {\em Numerical Linear Algebra}, vol.~50,
  SIAM, 1997.

\bibitem{van1935verteilungsfunktionen}
{\sc J.~van~der Corput}, {\em Verteilungsfunktionen}, in Proc. Akad. Amsterdam,
  vol.~38, 1935, p.~6.

\bibitem{weyl1916gleichverteilung}
{\sc H.~Weyl}, {\em {\"U}ber die {G}leichverteilung von {Z}ahlen mod. eins},
  Mathematische Annalen, 77 (1916), pp.~313--352.

\bibitem{wilson1970discrete}
{\sc M.~W. Wilson}, {\em Discrete least squares and quadrature formulas},
  Mathematics of Computation, 24 (1970), pp.~271--282.

\bibitem{wilson1970necessary}
{\sc M.~W. Wilson}, {\em Necessary and sufficient conditions for equidistant
  quadrature formula}, SIAM Journal on Numerical Analysis, 7 (1970),
  pp.~134--141.

\end{thebibliography}

\end{document}